\documentclass[11pt,reqno]{amsart}
\usepackage[normalem]{ulem}
\PassOptionsToPackage{table}{xcolor}
\usepackage{stmaryrd}
\expandafter\def\csname opt@stmaryrd.sty\endcsname
{only,shortleftarrow,shortrightarrow}
\usepackage{extpfeil}
\usepackage{amsmath,amssymb,amsthm}
\usepackage{aliascnt}
\usepackage{varioref}

\usepackage[backref=page]{hyperref}   

\renewcommand*{\backref}[1]{}          
\renewcommand*{\backrefalt}[4]{%
  \ifcase #1 %
    \relax                             
  \or
    Cited on page~#2.%
  \else
    Cited on pages~#2.%
  \fi}

\usepackage{thmtools}
\usepackage{microtype}
\usepackage[capitalize,nameinlink,noabbrev]{cleveref}
\hypersetup{colorlinks=true,urlcolor=blue,citecolor=blue,linkcolor=blue}
\usepackage{courier}
\usepackage{tikz}
\usepackage{tikz-cd}
\usetikzlibrary{calc,matrix,arrows,decorations.markings}
\usepackage{array}
\usepackage{color}
\usepackage{cancel}
\usepackage{enumerate}
\usepackage{nicefrac}
\usepackage{listings}
\usepackage{seqsplit}
\usepackage{longtable}
\usepackage{colonequals}
\usepackage{wasysym}
\usepackage{array, multirow}
\usepackage{mathtools}
\usepackage{quiver}

\newcommand{\repo}[2]{\href{https://github.com/fsaia/Unexpected-Isos/blob/main/#1}{#2}}
\newcommand{\repoDistinctProducts}[1]{\repo{Isogenies/DifferentProductLevel/#1}{\nolinkurl{#1}}}

\newcommand{\mathsout}[1]
{\bgroup\mathchoice
  {\sbox0{$\displaystyle{#1}$}%
    \usebox0\hspace{-\wd0}%
    \rule[0.5\ht0-0.5\dp0-.5pt]{\wd0}{1pt}}%
  {\sbox0{$\textstyle{#1}$}%
    \usebox0\hspace{-\wd0}%
    \rule[0.5\ht0-0.5\dp0-.5pt]{\wd0}{1pt}}%
  {\sbox0{$\scriptstyle{#1}$}%
    \usebox0\hspace{-\wd0}%
    \rule[0.5\ht0-0.5\dp0-.5pt]{\wd0}{1pt}}%
  {\sbox0{$\scriptscriptstyle{#1}$}%
    \usebox0\hspace{-\wd0}%
    \rule[0.5\ht0-0.5\dp0-.5pt]{\wd0}{1pt}}%
\egroup}

\usepackage[margin=1.05in]{geometry}

\newcommand{\Fp}{\mathbb{F}_p}
\newcommand{\Fpbar}{\overline{\mathbb{F}}_p}

\newcommand{\F}{\mathbb{F}}
\newcommand{\Q}{\mathbb{Q}}

\newcommand{\R}{\mathbb{R}}
\newcommand{\C}{\mathbb{C}}

\newcommand{\Z}{\mathbb{Z}}

\newcommand{\Zhat}{\widehat\Z}

\newcommand{\Aut}{{\rm Aut}}
\newcommand{\Gal}{{\rm Gal}}
\newcommand{\SL}{{\rm SL}}
\newcommand{\GL}{{\rm GL}}
\newcommand{\PGL}{{\rm PGL}}
\newcommand{\PSL}{{\rm PSL}}

\newcommand{\Frob}{{\rm Frob}}
\newcommand{\Ver}{{\rm Ver}}

\newcommand{\id}{\operatorname{id}}

\renewcommand{\P}{\mathbb{P}}

\newcommand{\End}{\operatorname{End}}

\newcommand{\Spec}{\operatorname{Spec}}

\newcommand{\Jac}{\operatorname{Jac}}

\newcommand{\new}{\operatorname{new}}

\DeclareFontFamily{U}{wncy}{}
    \DeclareFontShape{U}{wncy}{m}{n}{<->wncyr10}{}
    \DeclareSymbolFont{mcy}{U}{wncy}{m}{n}
    \DeclareMathSymbol{\Sh}{\mathord}{mcy}{"58}

\newtheorem{theorem}{Theorem}[section]

\newtheorem{lemma}[theorem]{Lemma}
\newtheorem{proposition}[theorem]{Proposition}
\newtheorem{corollary}[theorem]{Corollary}

\theoremstyle{definition}
\newtheorem{definition}[theorem]{Definition}
\newtheorem{example}[theorem]{Example}

\newtheorem{question}[theorem]{Question}

\newtheorem{remark}[theorem]{Remark}

\newtheorem{setup}[theorem]{Setup}
\crefname{setup}{Setup}{Setups}

\usepackage{xcolor}

\begin{document}
\title{Unexpected isogenies and isomorphisms between quotients of modular and Shimura curves}

\begin{abstract}
One can quotient a modular curve $X_0(N)$, or more generally a Shimura curve $X_0^D(N)$, by a subgroup of its group of Atkin--Lehner involutions to obtain a new curve. For squarefree levels $N$, we give a full classification of pairs of such quotients whose Jacobians are isogenous by relating the existence of such an isogeny to the existence of a newform of specified level with prescribed sign pattern. We also give a full classification of isomorphisms over $\mathbb{Q}$ between such quotients of the modular curves $X_0(N)$, and we work towards a complete classification of isomorphisms between such quotients of the curves $X_0^D(N)$. Part of this work involves describing the reduction type of an Atkin--Lehner quotient of $X_0^D(N)$ at a prime $p \parallel N$, and along the way we use this description to determine new defining equations for six Shimura curves quotients. 
\end{abstract}


\author{\sc Nikola Ad\v zaga}
\address{Nikola Ad\v zaga \\
University of Zagreb Faculty of Civil Engineering\\
Croatia}
\urladdr{https://www.grad.unizg.hr/en/nikola.adzaga}
\email{nikola.adzaga@grad.unizg.hr}

\author{\sc Oana P\u{a}durariu}
\address{Oana P\u{a}durariu \\
Max-Planck-Institut für Mathematik Bonn\\
Germany}
\urladdr{https://sites.google.com/view/oanapadurariu/home}
\email{opadurariu@mpim-bonn.mpg.de}

\author{\sc Frederick Saia}
\address{Frederick Saia \\
Lafayette College\\
USA}
\urladdr{https://fsaia.github.io/site/}
\email{saiaf@lafayette.edu}

\maketitle



\section{Introduction}\label{section: introduction}

Mazur's Program B \cite[p.~109]{Mazur_rat_pts} asks, given a number field $K$ and a subgroup $H \leq \GL_2(\Zhat)$, for a classification of all elliptic curves $E/K$ whose adelic Galois representation \[ \rho_E \colon \Gal(\overline{K}/K) \to \GL_2(\Zhat)\] has image conjugate to a subgroup of $H$. Equivalently, this program seeks a classification of the non-cuspidal rational points on modular curves over number fields. The program began with Mazur's determination of $X_1(N)(\Q)$ \cite{Mazur_X1} and of $X_0(N)(\Q)$ for prime $N$ \cite{Mazur_isogenies}. It has since grown to include Atkin--Lehner quotients of modular curves (e.g., \cite{Momose86, Galbraith02, BBBLMTV, BGX21, ACKP, AABCCKW, AM23, BDMTV, ADK26}) as well as Shimura curves and their Atkin--Lehner quotients (e.g., \cite{Clark03, BFGR06, RDVP14, NR15, Arai18, Arai23, PS23, LSSV24, PS25}), driven by interest in $\Q$-curves and in abelian surfaces with potential quaternionic multiplication. Such curves are often handled with sophisticated tools like the Chabauty--Coleman \cite{Coleman} and quadratic Chabauty \cite{BDMTV} methods. Sometimes a much simpler approach works: exhibit an isomorphism $X \to Y$, where $Y/K$ is a curve whose $K$-rational points are already known. Indeed, isomorphisms over $\Q$ between Atkin--Lehner quotients have been used in the literature to avoid redundant work, as the following examples show.

\begin{itemize}
    \item In \cite[Section 3.1]{ACKP}, the authors compute the $\mathbb{Q}$-rational points on the full Atkin--Lehner quotients $X_0(N)^*$ for all levels $N$ for which this quotient is hyperelliptic. Letting $X_0(N)^+ \coloneqq X_0(N)/\langle w_N \rangle$ denote the quotient by the Fricke involution, they note the isomorphisms
    \[ X_0(2p)^* \cong X_0(p)^+  \quad \text{ for } p \in \{67,73,103\} \quad \text{ and } \quad X_0(266)^* \cong X_0(133)^*, \]
    which reduce the rational points question for the former $3$ curves to ones handled in \cite{BBBLMTV} and \cite{BDMTV} and reduce that of $X_0(266)^*$ to the curve $X_0(133)^*$ which is also in their purview. 
    \item In \cite{PS23}, the authors determine the rational points on all Atkin--Lehner quotients of geometrically hyperelliptic Shimura curve $X_0^D(N)$ with $D>1$. They use the following isomorphisms to quotients of modular curves whose rational points have previously been computed (namely, in \cite{BGX21, BBBM, ACKP}):
    \[ X_0^{93}(1)^+ \cong X_0(93)^*, \quad X_0^{91}(1)^+ \cong X_0(91)^+, \quad \text{ and } \quad X_0^{85}(1)^+ \cong X_0(255)^*. \]
\end{itemize}

In this work, we are motivated by the goal of determining all isomorphisms between distinct Atkin--Lehner quotients of the Shimura curves $X_0^D(N)$. This includes the $D=1$ case, which we take to denote the compact modular curve $X_0^1(N) \coloneqq X_0(N)$ (rather than the curve $Y_0(N)$, as is also a reasonable convention). We note here that we consider general isomorphisms between these curves as algebraic curves over $\mathbb{Q}$, and in particular we \emph{are not} restricting to those that are modular maps (see \cref{remark: exceptional}).

Our first main result is a complete answer, in the squarefree level case, to a coarser problem; if $X$ and $Y$ are two Atkin--Lehner quotients which are isomorphic over $\mathbb{Q}$, or more generally which are isomorphic over $\overline{\mathbb{Q}}$, then it is necessary that their Jacobians $\Jac(X)$ and $\Jac(Y)$ are isogenous, and we determine all such isogenies for $X$ and $Y$ of squarefree levels. In the genus $1$ case, isomorphism classes of Jacobians of elliptic curves can be found in \cref{table: genus_1_modular_quotients} for the $D=1$ case and in \cite[Theorem 1.2]{PS25} in the $D>1$ case, so we restrict to genus $g \geq 2$ in the statements of our main theorems. 

\begin{theorem}\label{theorem: main_isogenies_theorem}
Let $(D_1, N_1, W_1)$ and $(D_2,N_2,W_2)$ be distinct triples such that 
\begin{itemize}
    \item $D_i$ is the discriminant of an indefinite quaternion algebra over $\Q$, 
    \item $N_i$ is a squarefree positive integer that is relatively prime to $D_i$, and 
    \item $W_i$ is a subgroup of the group of Atkin--Lehner involutions of $X_0^{D_i}(N_i)$
\end{itemize}
for $i \in \{1,2\}$. If there exists an isogeny 
\[ \Jac\left(X_0^{D_1}(N_1)/W_1\right) \sim \Jac\left(X_0^{D_2}(N_2)/W_2 \right) \]
over $\overline{\Q}$, then there exists such an isogeny over $\Q$, and there are only finitely many such isogenous pairs. Furthermore, we have $\max\{D_1N_1,D_2N_2\} \leq 4290$, and each such pair of isogenous Jacobians of quotients of genus $g \geq 2$ is recorded in \cref{table: isogenies_same_curve}, \cref{table: isogenies_same_product} or \cref{tab:isog_shimura_different_levels}.
\end{theorem}

While \cref{theorem: main_isogenies_theorem} does more in explicitly determining \emph{all} pairs of distinct Atkin--Lehner quotients of Shimura curves with isogenous Jacobians, we highlight that in particular it states the finitude of such pairs. It is not a-priori clear that one has finitude in squarefree level, and in fact finitude does \emph{not} hold in general level (see \cref{FH_isos}), but it should be expected from heuristics on dimensions of spaces of newforms with prescribed Atkin--Lehner actions coming from work of Martin \cite{Martin18}. These heuristics motivate a recent result of Padurariu--Park--Voight \cite{PPV26}, which proves explicit finitude of such spaces of dimension $0$ and which will be integral to our proof of \cref{theorem: main_isogenies_theorem}. The referenced result is restricted to squarefree levels, hence the restriction in \cref{theorem: main_isogenies_theorem}. The largest level $D_iN_i$ appearing for a Jacobian in an isogenous pair, as mentioned in \cref{theorem: main_isogenies_theorem}, is $4290$, which by \cite[Theorem 1.1]{PPV26} is the largest level $M$ for which there exists a sign pattern 
admitting a zero-dimensional space of newforms. 

Following \cref{theorem: main_isogenies_theorem}, we are able to reach a partial result concerning isomorphisms between Atkin--Lehner quotients over $\mathbb{Q}$. 

\begin{theorem}\label{theorem: main_isomorphisms_theorem}
We have the following isomorphisms over $\mathbb Q$ of
Atkin--Lehner quotients of Shimura curves of genus $2$:
\begingroup
\small
\allowdisplaybreaks
\begin{gather*}
X_0(26) \cong X_0^{26}(7)/\langle w_{182}\rangle \cong X_0^{26}(15)/\langle w_{6},w_{65}\rangle \\
X_0(30)/\langle w_{2}\rangle \cong X_0^{15}(22)/\langle w_{6},w_{55}\rangle \cong X_0^{15}(34)/\langle w_{2},w_{255}\rangle \\
X_0(93)/\langle w_{3},w_{31}\rangle \cong X_0(186)/\langle w_{2},w_{3},w_{31}\rangle \cong X_0^{93}(1)/\langle w_{93}\rangle \\
X_0(170)/\langle w_{2},w_{5},w_{17}\rangle \cong X_0(255)/\langle w_{3},w_{5},w_{17}\rangle \cong X_0^{85}(1)/\langle w_{85}\rangle \\
X_0(42)/\langle w_{3}\rangle \cong X_0^{14}(15)/\langle w_{6},w_{35}\rangle \cong X_0^{14}(39)/\langle w_{3},w_{182}\rangle
\end{gather*}
\endgroup

\begingroup
\renewcommand{\arraystretch}{1.35}

\begin{longtable}{c@{\quad}|@{\quad}c}
\(X_0(46)/\langle w_{46}\rangle \cong X_0(69)/\langle w_{69}\rangle\)
&
\(X_0(57)/\langle w_{57}\rangle \cong X_0(114)/\langle w_{3},w_{38}\rangle\)
\\
\(X_0(58)/\langle w_{29}\rangle \cong X_0^{58}(3)/\langle w_{174}\rangle\)
&
\(X_0(66)/\langle w_{11}\rangle \cong X_0^{6}(55)/\langle w_{11},w_{30}\rangle\)
\\
\(X_0(67)/\langle w_{67}\rangle \cong X_0(134)/\langle w_{2},w_{67}\rangle\)
&
\(X_0(70)/\langle w_{35}\rangle \cong X_0^{14}(15)/\langle w_{10},w_{21}\rangle\)
\\
\(X_0(73)/\langle w_{73}\rangle \cong X_0(146)/\langle w_{2},w_{73}\rangle\)
&
\(X_0(74)/\langle w_{74}\rangle \cong X_0(111)/\langle w_{111}\rangle\)
\\
\(X_0(78)/\langle w_{39}\rangle \cong X_0^{26}(3)/\langle w_{78}\rangle\)
&
\(X_0(78)/\langle w_{2},w_{13}\rangle \cong X_0^{39}(2)/\langle w_{3},w_{26}\rangle\)
\\
\(X_0(91)/\langle w_{91}\rangle \cong X_0^{91}(1)/\langle w_{91}\rangle\)
&
\(X_0(102)/\langle w_{3},w_{17}\rangle \cong X_0^{34}(3)/\langle w_{102}\rangle\)
\\
\(X_0(103)/\langle w_{103}\rangle \cong X_0(206)/\langle w_{2},w_{103}\rangle\)
&
\(X_0(115)/\langle w_{5},w_{23}\rangle \cong X_0(230)/\langle w_{2},w_{5},w_{23}\rangle\)
\\
\(X_0(133)/\langle w_{7},w_{19}\rangle \cong X_0(266)/\langle w_{2},w_{7},w_{19}\rangle\)
&
\(X_0^{6}(55)/\langle w_{2},w_{165}\rangle \cong X_0^{33}(2)/\langle w_{33}\rangle\)
\\
\(X_0^{6}(55)/\langle w_{6},w_{55}\rangle \cong X_0^{22}(15)/\langle w_{10},w_{11}\rangle\)
&
\(X_0^{6}(85)/\langle w_{6},w_{85}\rangle \cong X_0^{51}(2)/\langle w_{51}\rangle\)
\\
\(X_0^{10}(19)/\langle w_{190}\rangle \cong X_0^{38}(5)/\langle w_{5},w_{38}\rangle\)
&
\(X_0^{10}(21)/\langle w_{2},w_{3}\rangle \cong X_0^{35}(6)/\langle w_{15},w_{21}\rangle\)
\\
\(X_0^{10}(21)/\langle w_{10},w_{14}\rangle \cong X_0^{15}(14)/\langle w_{3},w_{7}\rangle\)
&
\(X_0^{10}(33)/\langle w_{2},w_{165}\rangle \cong X_0^{55}(2)/\langle w_{55}\rangle\)
\\
\(X_0^{10}(39)/\langle w_{5},w_{78}\rangle \cong X_0^{65}(2)/\langle w_{5},w_{13}\rangle\)
&
\(X_0^{10}(51)/\langle w_{3},w_{170}\rangle \cong X_0^{34}(15)/\langle w_{3},w_{10},w_{17}\rangle\)
\\
\(X_0^{10}(57)/\langle w_{2},w_{285}\rangle \cong X_0^{38}(5)/\langle w_{2},w_{19}\rangle\)
&
\(X_0^{10}(77)/\langle w_{10},w_{77}\rangle \cong X_0^{22}(5)/\langle w_{11}\rangle\)
\\
\(X_0^{14}(33)/\langle w_{2},w_{7},w_{11}\rangle \cong X_0^{462}(5)/\langle w_{2},w_{7},w_{165}\rangle\)
&
\(X_0^{14}(51)/\langle w_{2},w_{7},w_{51}\rangle \cong X_0^{14}(85)/\langle w_{2},w_{7},w_{85}\rangle\)
\\
\(X_0^{15}(2)/\langle w_{2}\rangle \cong X_0^{15}(14)/\langle w_{6},w_{7}\rangle\)
&
\(X_0^{21}(10)/\langle w_{7},w_{15}\rangle \cong X_0^{210}(11)/\langle w_{2},w_{7},w_{165}\rangle\)
\\
\(X_0^{21}(10)/\langle w_{14},w_{15}\rangle \cong X_0^{21}(22)/\langle w_{6},w_{154}\rangle\)
&
\(X_0^{22}(3)/\langle w_{2}\rangle \cong X_0^{33}(10)/\langle w_{15},w_{33}\rangle\)
\\
\(X_0^{22}(15)/\langle w_{5},w_{22}\rangle \cong X_0^{33}(2)/\langle w_{3}\rangle\)
&
\(X_0^{22}(13)/\langle w_{286}\rangle \cong X_0^{26}(11)/\langle w_{11},w_{26}\rangle\)
\\
\(X_0^{22}(23)/\langle w_{2},w_{253}\rangle \cong X_0^{46}(11)/\langle w_{2},w_{253}\rangle\)
&
\(X_0^{38}(5)/\langle w_{2},w_{95}\rangle \cong X_0^{38}(15)/\langle w_{2},w_{5},w_{57}\rangle\)
\\
\end{longtable}

\endgroup

of genus $3$:
\begingroup
\allowdisplaybreaks
\[
\renewcommand{\arraystretch}{1.35}
\begin{array}{c@{\quad}|@{\quad}c}
X_0(58)/\langle w_{58}\rangle \cong X_0(174)/\langle w_{2},w_{87}\rangle & X_0(78)/\langle w_{26}\rangle \cong X_0^{39}(2)/\langle w_{78}\rangle \\
X_0(97)/\langle w_{97}\rangle \cong X_0(194)/\langle w_{2},w_{97}\rangle & 
\end{array}
\] \endgroup

and of genus $4$:
\begingroup
\allowdisplaybreaks
\[
\renewcommand{\arraystretch}{1.35}
\begin{array}{c@{\quad}|@{\quad}c}
X_0(93)/\langle w_{93}\rangle \cong X_0(186)/\langle w_{3},w_{62}\rangle & X_0(170)/\langle w_{5},w_{34}\rangle \cong X_0(255)/\langle w_{5},w_{51}\rangle
\end{array}
\]
\endgroup
Furthermore, if $(D_1,N_1,W_1)$ and $(D_2,N_2,W_2)$ are distinct
triples as in the statement of
\cref{theorem: main_isogenies_theorem} and there exists an isomorphism
\[
X_0^{D_1}(N_1)/W_1\cong X_0^{D_2}(N_2)/W_2
\]
of curves of genus $g\geq 2$ over $\mathbb Q$ that is not among those listed above, then this pair of curves
must appear among the $151$ pairs listed in
\cref{table: remaining_isomorphism_candidates_same_product}
and \cref{table: remaining_isomorphism_candidates_different_products}.
\end{theorem}

The lack of a complete answer in \cref{theorem: main_isomorphisms_theorem} comes down to a lack of knowledge of equations for Atkin--Lehner quotients of Shimura curves; for each pair in \cref{table: remaining_isomorphism_candidates_same_product} and in \cref{table: remaining_isomorphism_candidates_different_products}, at least one curve is a Shimura curve quotient for which we lack an explicit defining equation over $\mathbb{Q}$. Hence, while none of the tests we implement exclude the possibility of isomorphism over $\mathbb{Q}$ for the pair, we lack a means for \emph{certifying} isomorphism over $\mathbb{Q}$. That said, \cref{theorem: main_isomorphisms_theorem} gives a satisfyingly thorough answer for modular curves. 

\begin{corollary}\label{theorem: main_isomorphisms_theorem_modular}
The only isomorphisms over $\Q$ between distinct Atkin--Lehner quotients of modular curves $X_0(N)$, with squarefree level $N$ and genus greater than $1$, are those listed in \cref{theorem: main_isomorphisms_theorem}. In particular, there are $10$ pairs of isomorphic genus $2$ curves, $2$ pairs of isomorphic genus $3$ curves, $2$ pairs of isomorphic genus $4$ curves, and no pairs of isomorphic curves having genus $g>4$.
\end{corollary}

One notices from this result that there exist no isomorphisms over $\Q$ between distinct Atkin--Lehner quotients of the same modular curve $X_0(N)$ in genus $g>1$ and in squarefree level $N$. 
\begin{corollary}
Let $N$ be a positive squarefree integer and $W, W'$ distinct Atkin--Lehner subgroups of $\Aut \, X_0(N)$. There are no isomorphisms over $\Q$ of the form $X_0(N)/W \cong X_0(N)/W'$ between curves of genus $g > 1$.
\end{corollary}
In genus $1$, there are exactly two such isomorphisms of quotients of modular curves in squarefree level; we have
\begin{align*}
    X_0(30)/\langle w_6 \rangle &\cong X_0(30)/\langle w_3, w_{10} \rangle \quad \text{ and } \\
     X_0(42)/\langle w_{14} \rangle &\cong X_0(42)/\langle w_6, w_{7} \rangle,
\end{align*} 
with the former pair having Cremona reference $15$a$3$ and the latter pair having reference $21$a$4$ (see \cref{table: genus_1_modular_quotients}). 

\cref{theorem: main_isomorphisms_theorem_modular} improves upon work of Padurariu in \cite[Section 4.1]{Oana_thesis}, wherein the author gives a classification of isomorphisms within some restricted families of Atkin--Lehner quotients of modular curves. Namely, she handles the case of $X_0(p)^+ \cong X_0(2p)^+$ for prime level $p$ and that of $X_0(pq)^* \cong X_0(q)^*$ for primes $p$ and $q$. While these cases, and other restricted families, can be reduced to finite problems just using genus considerations, it is important to note that this is \emph{not} possible for the full generality of main results appearing above. We comment further on this point in \cref{subsection:genus_formulae}. 

One can extend our questions and test our methods in various directions, such as removing the squarefree condition, or working with more general congruence subgroups. For a more detailed discussion, see \Cref{section: Zukunftsmusik}.

We now briefly describe the theory we developed and applied in order to prove our main results.
To prove \cref{theorem: main_isogenies_theorem}, we use Ribet's isogeny to show that the vanishing of certain spaces of newforms with prescribed sign patterns is necessary for an isogeny between Jacobians of Atkin--Lehner quotients. Together with the classification of such vanishing spaces by Padurariu--Park--Voight \cite{PPV26}, these steps bound levels admitting isogenies. Furthermore, we use descriptions coming from the moduli interpretation, e.g., point counts of Shimura curve quotients over finite fields, studied and implemented in \cite{MPSS25}. In the other direction, while in some cases we exhibit isogenies by computing Jacobian decompositions, in many cases we do without full decompositon computations by comparing the $p$-old part of the Jacobian with the way sign patterns extend from level $DN$ to level $DNp$ (\cref{lem:p-old-criterion}, \cref{lem:change_in_D_criterion}).

A main obstacle in transitioning to \cref{theorem: main_isomorphisms_theorem} is that equations for (Atkin--Lehner quotients of) Shimura curves are largely unknown for $D>1$. We therefore use information, particularly local invariants, that come from the moduli interpretation alone. Work of Ogg \cite{Ogg83} provides information on real points, while at primes $p\mid D$ we have reduction type descriptions as implemented in \cite{PS25} \`{a} la {\v{C}}erednik--Drinfeld \cite{Cerednik, Drinfeld}. A new tool we implement, generalizing the $D=1$ case that was worked out in \cite{ADK26}, is the description of a semi-stable model of the quotients $X_0^D(N)/W$ for primes $p \parallel \,N$ (in \Cref{thm: semistable_XD_quotient}) and its minimal desingularization (in \Cref{cor:minimal_regular_AL_quot}) \`{a} la Deligne--Rapoport \cite{DeligneRapoport}. For example, we prove the following as a consequence.

\begin{theorem}\label{thm: stable}
Let $N$ be a positive integer and $p$ be a prime that exactly divides $N$, and $D$ be the discriminant of an indefinite quaternion algebra over $\Q$, such that $D$ and $N$ are coprime. Assume that $X_0^D(N)/W$ has genus at least $2$, where $W$ is a subgroup of the group of Atkin--Lehner involutions. Then the curve $X_0^D(N)/W$ has stable reduction over $\Z_p$, and hence its minimal regular model is semi-stable.
\end{theorem}

As additional applications of our descriptions of these models, we give an equivalent criterion for the existence of $\Q_p$-points on $X_0^D(N)/W$ for primes $p\parallel N$ (in \Cref{cor: existence_of_Qp_points}) and we determine $6$ new equations for Atkin--Lehner quotients of $X_0^D(N)$ over $\Q$ (see \cref{prop : bielliptic_improvement} and \cref{prop : non-elliptic_improvement}). 

The outline of this paper is as follows. In \cref{section: background}, we discuss relevant background on modular forms and modular abelian varieties, modular and Shimura curves, and Atkin--Lehner involutions thereof. We proceed to discussing isogenous Jacobians of Atkin--Lehner quotients in \cref{section: isogenies}, recalling the above results of Padurariu--Park--Voight referenced above in \cref{section: PPV}. This section contains the proof of \cref{theorem: main_isogenies_theorem}, broken into several cases in \cref{subsection: isogenies_same_curve}, \cref{subsection: isogenies_same_level}, and \cref{subsection: isogenies_different_level}. After handling reduction types of Atkin--Lehner quotients at primes dividing the level in \cref{section: reductions_at_N},  we prove \cref{theorem: main_isomorphisms_theorem} in \cref{section: isomorphisms}. In \cref{section: Zukunftsmusik}, we briefly comment on possible extensions of this work, namely to general level or to other families of Shimura curves. On the former topic, we remark that while there \emph{are infinitely many} isomorphisms between distinct Atkin--Lehner quotients over $\mathbb{Q}$, our results from this work should still prove useful towards classifications of isogenies and isomorphisms in non-squarefree level. 

All computations described in this paper were performed in Magma \cite{magma}. All
related code, as well as computed data, can be found in the Github repository \cite{APS_github}.

\section*{Acknowledgments} 
We are pleased to thank Sachi Hashimoto, Sun Woo Park, and John Voight for helpful conversations. This project was completed in part while N.A. was a visiting scientist at the Max-Planck-Institut f\"{u}r Mathematik in Bonn, and he is grateful to the institute for its support and hospitality.

N.A.\@ was supported by an institutional project of the University of Zagreb Faculty of Civil Engineering (GEM -- Graphs, elliptic and modular curves), under the programme line for enhancing scientific excellence, financed through the National Recovery and Resilience Plan 2021–2026. F.S.\@ acknowledges support for this research provided by an AMS--Simons Travel Grant. 

\section*{AI Statement} No AI was used in developing the mathematical arguments appearing in this work. Generative AI was employed to draft some of the initial code, but this code was fully rewritten and verified by the authors. AI was also used to proofread a draft of the background section (\cref{section: background}) for errors. 
After using these tools, the authors reviewed and edited the content as needed and take full responsibility for the content of this manuscript.


\section{Background}\label{section: background}

\subsection{Modular and Shimura curves}\label{subsection: modular and Shimura curves}

Let $D \geq 1$ be the discriminant of an indefinite quaternion algebra $B$ over $\mathbb{Q}$, i.e., $D$ is the product of an even number of distinct prime numbers, and let $N$ be a positive integer which is relatively prime to $D$. We let $X_0^D(N)$ denote the canonical model over $\mathbb{Q}$, as given by Shimura \cite[Main Theorem 1]{Sh67}, for the coarse moduli scheme parameterizing pairs $(A, \iota_N)$ where $A$ is an abelian surface and $\iota_N: O_N \hookrightarrow \End(A)$ is an embedding of an Eichler order $O_N$ of $B$. Equivalently (up to isomorphism of the coarse moduli scheme), this curve parameterizes triples $(A,\iota,Q)$ where $A$ is an abelian surface, $\iota: \mathcal{O} \hookrightarrow \End(A)$ is an embedding of a maximal order $\mathcal{O}$ of $B$, and $Q \leq A[N]$ is a cyclic $\mathcal{O}$-submodule which as a group is isomorphic to $\left(\Z/N\Z\right)^2$.\footnote{Here we should also remark that part of the datum is a principal polarization on $A$ that is compatible with the QM-structure $\iota$. By \cite{Sh67} there is a unique such polarization up to a choice that does not change the isomorphism class over $\Q$ of the coarse moduli scheme, permitting us to suppress the polarization.} For background on quaternion algebras, we refer the reader to \cite{Voight}.

If $D = 1$, we take the convention that this curve is the compactified modular curve $X_0^1(N) \coloneqq X_0(N)$ (rather than the non-compact curve $Y_0(N)$, which is also a common convention). The curve $X_0(N)$, up to isomorphism, can also be described as the coarse moduli scheme over $\mathbb{Q}$ parameterizing pairs $(E,C)$ where $E$ is an elliptic curve and $C \leq E[N]$ is a cyclic subgroup of order $N$.

\subsection{Modular forms, old and new}

Let $N$ be a positive integer, and consider the level $N$ congruence subgroup
\[ \Gamma_0(N) \coloneqq \left\{ \begin{bmatrix} a & b \\ c & d \end{bmatrix} \in \SL_2(\Z) : \begin{bmatrix} a & b \\ c & d \end{bmatrix} \equiv \begin{bmatrix} * & * \\ 0 & * \end{bmatrix} \pmod{N} \right \} \]
For $k$ a positive integer, in this work we denote by 
\[ M_k(N) \coloneqq M_k(\Gamma_0(N)) \, \text{ and } \,
S_k(N) \coloneqq S_k(\Gamma_0(N)) \]
the spaces of modular forms of weight $k$ and of cusp forms of weight $k$ with respect to $\Gamma_0(N)$, respectively (see \cite[Definition 1.2.3]{DS16}). 

Let $M \mid N$ and let $e$ be a divisor of $N/M$. We have an embedding
\begin{align*}
    \iota_e : S_k(M) &\hookrightarrow S_k(N) \\
    f(z) &\mapsto  f(ez).
\end{align*}
For a divisor $n \mid N$, the space of $n$-old cusp forms is the subspace of $S_k(N)$ spanned by the images of these embeddings for $M \mid N$ not divisible by $n$:
\[ S_k^{n-\text{old}}(N) \coloneqq \left \langle \bigcup_{\substack{M \mid N \\ n \nmid M}} \bigcup_{e \mid N/M} \text{Image}(\iota_e) \right \rangle. \]
The space $S_k^{n-\text{new}}(N)$ of new cusp forms is the orthogonal complement of $S_k^{n-\text{old}}(N)$ in $S_k(N)$ with respect to the Petersson inner product \cite[Definition 5.6.1]{DS16}. In particular, we set
\[ S_k^\text{old}(N) \coloneqq S_k^{N-\text{old}}(N) \,  \text{ and } \, 
S_k^\text{new}(N) \coloneqq S_k^{N-\text{new}}(N), \]
recovering the usual old/new decomposition at level $N$. The spaces $S_k^\text{old}(N)$ and $S_k^\text{new}(N)$ are stable under the action of the Hecke operators $T_n$ and $\langle n \rangle$ on $M_k(N)$ for each $n \in \mathbb{Z}^+$ \cite[Proposition 5.6.2]{DS16}. We refer to \cite[5.3]{DS16} for the definitions of these operators. 

\subsection{Newforms and modular abelian varieties}

By a \textbf{newform of level $N$}, we mean in this work a normalized simultaneous eigenform $f \in S_2(N)^\text{new}$ for the Hecke operators $T_n$ and $\langle n \rangle$ for all $n \in \mathbb{Z}^+$ (see \cite[Definition 5.8.1]{DS16}). In particular, we restrict from here on to weight $k=2$. To each newform $f$ of level $N$, we have an associated abelian variety $A_f$ over $\mathbb{Q}$ as given by Shimura \cite{Shi73} (see also \cite[Definition 6.63]{DS16}). We refer to any abelian variety arising in this way as a \textbf{modular abelian variety}.

Central to many arguments in this work, the Jacobian $\Jac(X_0(N))$ of $X_0(N)$ is isogenous to a direct sum of modular abelian varieties: letting $M_f$ denote the level of a newform $f$ and letting $d(n)$ denote the number of divisors of a positive integer $n$, we have (\cite[Theorem 6.6.6 and p. 250]{DS16})
\[ \Jac(X_0(N)) \sim \bigoplus_{\substack{\text{newforms } f \\ \text{of level $M_f \mid N$}}} A_f^{d(N/M_f)}, \]
where the sum is over a set of representatives $f$ for each orbit of newforms at level $M_f \mid N$ under the action of automorphisms of $\C$.

Concretely, for a newform $f$ of level $M \mid N$, there are $d(N/M)$ ways to embed $f$ into $S_2(N)$, and the image of $f$ under each embedding contributes an isogeny factor of $A_f$ in $\Jac(X_0(N))$. We collect some further results on these isogeny factors, due to Ribet, which will be crucial in what follows.

\begin{proposition}\label{Proposition: Ribet_modular_abelian_varieties}
Let $N$ be a squarefree positive integer and let $f \in S_2(N)$ be a newform of level $N$.
\begin{enumerate}
    \item The modular abelian variety $A_f$ is simple over $\overline{\mathbb{Q}}$. 
    \item If $g$ is a newform of squarefree level such that $A_f$ and $A_g$ are isogenous over $\overline{\mathbb{Q}}$, then $f$ and $g$ are in the same Galois orbit,  i.e., $A_f$ and $A_g$ are isogenous over $\mathbb{Q}$. 
\end{enumerate}
\end{proposition}
\begin{proof}
Part (1) follows from \cite[Corollary 1.4 and \S 3]{Ribet75}; the assumption that $N$ is squarefree is required to have that $\Jac(X_0(N))$ is semi-stable over $\mathbb{Q}$. Part (2) follows from \cite[Theorem (3.9 bis), Corollary (3.10), and Theorem (4.7)]{Ribet80}. 
\end{proof}

Further on, we always choose a single representative of each Galois orbit.

\subsection{Atkin--Lehner involutions, Ribet's isogeny, and compatibility}

For $\mathcal{O}_N$ an Eichler order of squarefree level $N$ in the quaternion algebra $B$ of discriminant $D$, we consider the following subgroup of the normalizer of $\mathcal{O}_N$:
\[N_{{B}^\times_{> 0}}(\mathcal{O}_N) \coloneqq \{u \in {B}^\times : \text{nrd}(u) > 0 \text{ and } u^{-1} \mathcal{O}_N u = \mathcal{O}_N \},\] 
where $\text{nrd}$ is the reduced norm. The group of \textbf{Atkin--Lehner involutions} on the curve $X_0^D(N)$ is the following subgroup of the automorphism group\footnote{In general in this work, for $X$ a variety over a base field $F$ we let $\Aut(X)$ denote the group of automorphisms of $X$ over $F$ (as opposed to those over $\overline{F}$).} over $\mathbb{Q}$:
\[ W_0(D,N) \coloneqq N_{{B}^\times_{> 0}}(\mathcal{O}_N) / \left(\mathbb{Q}^\times \mathcal{O}_N^\times \cap B^\times_{>0}\right) \leq \text{Aut}(X_0^D(N)). \]
Letting $\omega(n)$ denote the number of distinct prime divisors of a positive $n$, the Atkin--Lehner group of $X_0^D(N)$ contains a non-trivial involution $w_m$ for each Hall divisor $1 < m \parallel DN$ (i.e., $\gcd(m, DN/m)=1$) and moreover has the structure of a finite abelian $2$-group of the form
\[
W_0(D,N) = \{w_m : m \parallel DN\} \cong \left(\Z/2\Z\right)^{\omega(DN)}.
\]
For a subgroup $W \leq W_0(D,N)$, we refer to the quotient $X_0^D(N)/W$ as an \textbf{Atkin--Lehner quotient} of $X_0^D(N)$. We also refer to the elements of the subgroup of $\Aut\left(X_0^D(N)/W\right)$ induced by the elements of $W_0(D,N)$ as Atkin--Lehner involutions.

Based on the second moduli interpretation we gave for the curve $X_0^D(N)$, we can describe the action of an involution $w_{m}$ as follows in direct analogy to the classical modular curve case (see \cite[p. 83]{Elkies04}): for a point $P$ induced by a triple $(A,\iota,Q)$, we have two natural decompositions of the QM-cyclic isogeny
\[ \varphi : A \to A/Q \]
as 
\[ \varphi = \varphi_{N/m} \circ \varphi_m = \psi_m \circ \psi_{N/m}, \]
where $\varphi_m$ and $\psi_m$ are each isogenies of degree $m^2$. We let $Q' = \ker( \psi_{N/m} \circ (\varphi_m)^\vee ) \leq A/Q[m]$, and we note from the assumed properties on $Q$ that this subgroup $Q'$ is QM-stable and cyclic as an $\mathcal{O}$-module with respect to the natural QM-structure $\overline{\iota}$ induced by $\iota$ on $A/Q[m]$. The triple $(A/Q[m],\overline{\iota},Q')$ then induces a point on $X_0^D(N)$ which we defined to be $w_m(P)$. There is an equivalent description $w_m(A, Q) = (A/Q[m], Q+A[m])/Q[m])$ in which we suppress the QM-structure.

If $[A_f]$ is an isogeny factor of $\Jac(X_0(N))$, corresponding to a newform $f$, then as $W_0(1,N)$ acts on this Jacobian we can consider the action on $[A_f]$. We encode the information of this action in a tuple $\overline{\epsilon} = (\epsilon_1, \ldots, \epsilon_{\omega(N)})$ as follows: let $N = p_1 \cdots p_{\omega(N)}$ be the prime factorization of $N$, and let $\epsilon_i \in \{\pm 1\}$ be the eigenvalue of the action of $w_{p_i}$ on $f$. We call the tuple $\overline{\epsilon}$ the \textbf{sign pattern} associated to this isogeny factor.\footnote{Note that this sign pattern is only well defined up to the choice of ordering of the primes in the prime factorization of $N$. We will make clear our choice of ordering when appropriate.}

\begin{remark}\label{remark: sign_patterns_level_raising}
Note that if $A \sim A_f$ for a newform $f$ of level $M_f \mid N$, it \emph{does not} make sense to refer to this as \emph{the} sign pattern attached to $A_f$ or to $f$ \emph{unless} $M_f = N$; if $M_f$ is a proper divisor of $N$, then there are multiple ways to embed $f$ into $S_2(N)$. These embeddings give rise to distinct isogeny factors represented by $A_f$, and these distinct factors receive distinct actions by a fixed $w_m \in W_0(1,N)$. More specifically, the action of $w_m$ for each $m \mid M_f$ is common, determined by $f$, while the 
\[ d\left(N/M_f\right) = 2^{\omega(N/M_f)} \]
different embeddings of $f$ from $S_2(M_f)$ to $S_2(N)$ yield all possible choices of signs for the involutions $w_m$ with $m \mid (N/M_f)$. 
\end{remark}

The following result provides that the Jacobian of any Atkin--Lehner quotient of the Shimura curve $X_0^D(N)$ is isogenous to a specified sub-quotient of the Jacobian of the modular curve $X_0(DN)$. It was originally proven with level restrictions by Ribet \cite[Th\'eor\`{e}me 2]{Ribet80b}, which include the squarefree case and thus are suitable to our needs in the present work, though it can be extended to general level by results of Hijikata--Pizer--Shemanske (alternatively, see \cite[Corollary 3.3.2]{Roberts}). The claim regarding the Atkin--Lehner actions follows from \cite[\S 1.7]{BD96}.

\begin{theorem}\cite{Ribet80b, BD96, HPS89a, HPS89b}\label{Ribet_isog}
Let $D$ be the discriminant of an indefinite quaternion algebra over $\Q$ and let $N$ be a positive integer coprime to $D$ such that $X_0(DN)$ has positive genus. There is an isogeny defined over $\Q$
\[ \Psi \colon \Jac(X_0(DN))^{D-\new} \longrightarrow \Jac(X_0^D(N)) \]
such that, for each Atkin--Lehner involution $w_m(D,N) \in W_0(D,N)$, we have
\begin{equation}\label{eqn: AL_action_Ribet}
    \Psi^*(w_m(D,N)) = (-1)^{\omega(\gcd(D,m))}w_m(1,DN),
\end{equation} 
where $w_m(1,DN) \in W_0(1,DN)$ denotes the Atkin--Lehner involution corresponding to $m$ on the modular curve $X_0(DN)$.
\end{theorem}

It follows from \cref{Ribet_isog} that each simple isogeny factor of $\Jac(X_0^D(N))$ is represented by a modular abelian variety $A_f$ of a newform $f$ of some level $M_f \mid DN$. We can also consider the action of $W_0(D,N)$ on isogeny factors of $\Jac(X_0^D(N))$ as we did in the $D=1$ case above. The following definitions will help to simplify our discussion of those isogeny factors which are fixed by the action of specified Atkin--Lehner subgroups. 

\begin{definition}\label{definition: compatibility}
Let $D \geq 1$ be the discriminant of an indefinite quaternion algebra over $\Q$, let $N$ be a squarefree positive integer coprime to $D$, and let 
\[ \overline{\epsilon} = (\epsilon_1,\ldots, \epsilon_{\omega(DN)}) \in \{\pm 1\}^{\omega(DN)} \]
be a sign pattern in level $DN$. 

\begin{itemize}
    \item For $m \parallel DN$, we say that $\overline{\epsilon}$ is \textbf{compatible} with the action of the involution $w_m \in W_0(D,N)$ if 
\[ \prod_{p \mid m} \epsilon(p) = (-1)^{\omega(\gcd(D,m))}. \]
    
    \item For $W \leq W_0(D,N)$ an Atkin--Lehner subgroup, we say that $\overline{\epsilon}$ is \textbf{compatible} with the action of $W$ if $\overline{\epsilon}$ is compatible with the action of $w$ for every $w \in W.$
    
    \item For $M$ a divisor of $DN$ and $\overline{\epsilon}_M$ a sign pattern in level $M$, we say that $\overline{\epsilon}_M$ is \textbf{compatible} with the action of $w_m \in W_0(D,N)$ (respectively, with that of an Atkin--Lehner subgroup $W \leq W_0(D,N)$) if there exists a choice of sign at each prime $q \mid DN/M$ so that the resulting sign pattern $\overline{\epsilon}$ extending $\overline{\epsilon}_M$ is compatible with $w_m$ (respectively, with $W$).
    
    \item For $f$ a newform of level $M \mid DN$
with sign pattern $\overline{\epsilon}$, we say that $f$ and the corresponding modular abelian variety $A_f$ are \textbf{compatible} with the action of an involution $w_m \in W_0(D,N)$, or with that of a subgroup $W \leq W_0(D,N)$, if $\overline{\epsilon}$ is. 
\end{itemize}

\end{definition}

To spell out the worth of this terminology: by \cref{Ribet_isog} we have for $f$ as in the definition that $A_f$ represents an isogeny factor of $\Jac(X_0^D(N)/W)$ if and only if $f$ has level $M_f$ with $D \mid M_f$ and is compatible with the action of $W$. In this case, there may be multiple isogeny factors represented by $A_f$ if $M_f$ is a proper divisor of $DN$; the number of such factors is equal to the number of forms in $S_2(DN)$ with associated modular abelian variety $A_f$ (as mentioned in \cref{remark: sign_patterns_level_raising}) whose associated sign patterns in level $DN$ are compatible with the action of $W$, and so is at least $1$ and is at most $d(DN/M_f)$.

\begin{definition}
Let $D$ be the discriminant of an indefinite quaternion algebra over $\mathbb{Q}$, let $N$ be a squarefree positive integer coprime to $D$, and let $S \subseteq S_2(DN)$. 
\begin{itemize}
    \item For $\overline{\epsilon}$ a sign pattern in level $DN$, we let $S^{\overline{\epsilon}}$ denote the subspace spanned by the forms in $S$ with associated sign pattern $\overline{\epsilon}$. 
    \item For $w_m \in W_0(D,N)$ (respectively, for a subgroup $W \leq W_0(D,N)$) we let $S^{w_m}$ (respectively, $S^W$) denote the subspace spanned by the forms in $S$ which are compatible with the action of $w_m$ (respectively, of $W$). In other words, 
    \[ S^{W} \coloneqq \bigoplus_{\overline{\epsilon}}  S^{\overline{\epsilon}}, \]
    where the direct sum is over sign patterns $\overline{\epsilon}$ of level $DN$ which are compatible with the action of $W$.
\end{itemize}
In particular, for $n \mid DN$ and $\overline{\epsilon}$ a sign pattern in level $DN$ we let 
\[ S_2(DN)^{n-\text{new},\overline{\epsilon}} \coloneqq (S_2(DN)^{n-\text{new}})^{\overline{\epsilon}}, \]
denote the space spanned by those forms $f$ of level $M_f$ divisible by $n$ that have sign pattern $\overline{\epsilon}$, and similarly with $\overline{\epsilon}$ replaced by an involution $w_m \in W_0(D,N)$ or by a subgroup $W \leq W_0(D,N)$.
\end{definition}

\subsection{Genus formulae}\label{subsection:genus_formulae}
 Let \(D>1\) be the discriminant of an indefinite quaternion algebra over \(\Q\), and let \(N\) be squarefree and coprime to \(D\). Put
\(
\psi(N)=N\prod \limits_{\ell\mid N}\left(1+\frac{1}{\ell}\right).
\)
By \cite[Prop.~1.40]{ShimuraBook}, we have
\[
g(X_0^D(N))
=
1+\frac{\varphi(D)\psi(N)}{12}
-\frac{e_2(D,N)}{4}
-\frac{e_3(D,N)}{3},
\]
where
\[
e_2(D,N)=
\prod_{\ell\mid D}\left(1-\left(\frac{-4}{\ell}\right)\right)
\prod_{\ell\mid N}\left(1+\left(\frac{-4}{\ell}\right)\right),
\]
and
\[
e_3(D,N)=
\prod_{\ell\mid D}\left(1-\left(\frac{-3}{\ell}\right)\right)
\prod_{\ell\mid N}\left(1+\left(\frac{-3}{\ell}\right)\right).
\]
Here \(\left(\frac{\cdot}{\ell}\right)\) denotes the Kronecker symbol. In the classical \(D=1\) case, by \cite[Theorem 3.1.1]{DS16}, the genus formula has the same shape, but with the additional cusp contribution:
\[
g(X_0(N))
=
1+\frac{\psi(N)}{12}
-\frac{e_2(1,N)}{4}
-\frac{e_3(1,N)}{3}
-\frac{c(N)}{2},
\]
where
\( c(N)=\sum\limits_{d\mid N}\varphi(\gcd(d,N/d)) \).
In particular, if \(N\) is squarefree, then \(c(N)=2^{\omega(N)}\).

We will also use the Riemann--Hurwitz formula for Atkin--Lehner quotients. Let \(X \colonequals X_0^D(N)\), and let \(W\leq W_0(D,N)\). For a Hall divisor \(d\parallel DN\), let \(\nu^D(N,d)\) denote the number of fixed points of \(w_d\) on \(X_0^D(N)\) when \(D>1\). In the classical case \(D=1\), we write \(\nu(N, d) = \nu^1(N, d)\) for the number of fixed points of \(w_d\) on \(X_0(N)\).
The Riemann--Hurwitz formula gives
\[
g(X_0^D(N)/W)
=
1+\frac{g(X_0^D(N))-1}{|W|}
-\frac{1}{2|W|}
\sum_{\substack{d\parallel DN\\ 1 \neq w_d\in W}}
\nu^D(N,d).
\]

These formulae give a first, but rather crude, way of comparing different Atkin--Lehner quotients. For instance, if two quotients \(X_0^D(N)/W_1\) and \(X_0(DN)/W_2\) have the same order of Atkin--Lehner subgroup (\( |W_1|=|W_2| \)), then the equality of their genera gives a relationship between the  fixed-point sums
\[
\sum_{\substack{d\parallel DN\\ 1 \neq w_d\in W_1}}
\nu^D(N,d) \quad \text{ and} \quad
\sum_{\substack{d\parallel DN\\ 1 \neq w_d\in W_2}}
\nu(DN,d).
\]

\noindent The number of fixed points \( \nu^D(N, d)\) can be expressed by class numbers and Kronecker symbols, e.g., \( \displaystyle \nu(N, d) = (h(-d)+h(-4d))\prod_{p\mid N/d} \left(1+\left(\frac{-d}{p}\right)\right)\) for odd $N$ by \cite{Klu77}. One possible strategy is to try to compare the genus formulae using general upper and lower bounds for the class numbers. In some restricted families this gives useful inequalities. However, in the general situation we consider, the known bounds introduce error terms that can be too large compared with the main terms one is trying to separate. Thus, the bounds obtained in this way are not strong enough to reduce the general problem of determining all isogenies or isomorphisms between Atkin--Lehner quotients of Shimura curves to a finite computational box. 

Moreover, it is \emph{not possible} to reduce to a finite problem solely through genus considerations. For example, for \( m \equiv 1 \pmod{4}\), the number of fixed points of the involution $w_m$ on $X_0^D(N)$ with $m \parallel DN$ is
\[ \nu^D(N, m) = h(-4m)\prod_{p\mid D} \left(1- \left(\frac{-4m}{p}\right)\right)
\prod_{p\mid N} \left(1+ \left(\frac{-4m}{p}\right)\right).\]
One can then easily construct infinitely many squarefree pairs \( (D, N) \) such that \( g(X_0^D(N)/W_1) = g(X_0^D(N)/W_2) \) for two distinct Atkin--Lehner subgroups \( W_1, W_2\). 

\begin{example}
For every prime number $p \not \in \{5,7,13\}$, we have that the involutions $w_5$ and $w_{13}$ on the Shimura curve $X_0^{5\cdot 7 \cdot 13 \cdot p}(1)$ have no fixed points because 
\[ \left( \frac{-20}{7} \right) = \left( \frac{-52}{7} \right) = 1. \]
We then have 
\[ g\left( X_0^{5\cdot 7 \cdot 13 \cdot p}(1)/\langle w_5 \rangle \right) = g\left(X_0^{5\cdot 7 \cdot 13 \cdot p}(1) /\langle w_{13} \rangle \right). \]
\end{example}

\begin{example}
With the same idea, we see that in fact there are Shimura curves $X_0^D(N)$ with arbitrarily many Atkin--Lehner quotients of the same genus. Given an odd positive integer $r$, let $p_i$ be the $i^\text{th}$ smallest prime number satisfying $p_i \equiv 1 \pmod{4}$ for each $1 \leq i \leq r$. Let $q$ be a prime which is split in $\Q(\sqrt{-p_i})$ for all $1 \leq i \leq r$, and set 
\[ D \coloneqq p_1 \cdots p_r \cdot q_1. \]
Then the curves $X_0^D(1)/\langle w_{p_i} \rangle$ have the same genus for all $1 \leq i \leq r$. 
\end{example}

In the case of quotients of curves of different levels, we similarly see that genus considerations alone cannot reduce the general problem to consideration of finitely many curves. For example, there are infinitely many genera $g$ such that there exist distinct positive integers $N_1, N_2$ with 
\[ g = g(X_0(N_1)) = g(X_0(N_2)), \]
and moreover for each positive integer $r$ there exists a $g$ which is the genus of at least $r$ distinct modular curves $X_0(N)$ by \cite{CWZ00}. This all goes to show that genus formulae alone cannot reduce the problem of classifying isogenous Jacobians of Atkin--Lehner quotients to finitely many cases. Still, checking genus incidences serves as a useful coarse filter for the existence of isogenies after we narrow to a large, finite number of cases by means discussed in the next section.


\section{Isogenous Atkin--Lehner quotients}\label{section: isogenies}

In this section, we determine all tuples $(D_1,N_1,W_1, D_2,N_2,W_2)$ such that 
\begin{enumerate}
    \item $D_i \geq 1$ is the discriminant of an indefinite quaternion algebra over $\mathbb{Q}$,
    \item $N_i$ is a squarefree positive integer coprime to $D_i$, and
    \item $W_i \leq W_0(D_i,N_i)$ is an Atkin--Lehner subgroup
\end{enumerate}
for $i=1,2$, and such that there exists an isogeny
\[ \Jac\left(X_0^{D_1}(N_1)/W_1\right) \sim \Jac\left(X_0^{D_2}(N_2)/W_2\right) \]
over $\mathbb{Q}$. Equivalently, by \cref{Proposition: Ribet_modular_abelian_varieties}, we determine all such isogenous Jacobians over $\overline{\mathbb{Q}}$. Given this equivalence, we will refrain from further remarking on the field of definition of isogenies in this section. 

After a short aside on low-genus quotients in \cref{subsection: low_genus}, we will focus on Atkin--Lehner quotients of genus $g \geq 2$ and our determination will proceed in three cases:
\begin{itemize}
    \item the case of $(D_1,N_1) = (D_2,N_2)$ (that is, of quotients of the same Shimura curve) in \cref{subsection: isogenies_same_curve},
    \item the case of $D_1N_1 = D_2N_2$ with $D_1 \neq D_2$ (that of different Shimura curves of the same product level $DN$) in \cref{subsection: isogenies_same_level}, and
    \item the case of $D_1N_1 \neq D_2N_2$ (that of different product levels) in \cref{subsection: isogenies_different_level}.
\end{itemize}
Each case will see slightly different argumentation, and so we will state tools in the relevant sections. That said, before proceeding, we present here a result that is applicable in each of the first two cases. 

\begin{proposition}\label{prop: AL_restrictions}
Let $D_1$ and $D_2$ be discriminants of indefinite quaternion algebras over $\Q$ and let $N_1$ and $N_2$ be positive integers with $\gcd(D_i,N_i) = 1$ for $i \in \{1,2\}$ and with $D_1N_1 = D_2N_2$. Assume that there exists a newform of level $D_1N_1$ in every sign pattern. If there exist subgroups $W_1 \le W_0(D_1,N_1)$ and $W_2 \le W_0(D_2,N_2)$ such that 
\[ \Jac\left(X_0^{D_1}(N_1)/W_1\right) \sim \Jac\left(X_0^{D_2}(N_2)/W_2\right), \] 
then the following statements hold:
\begin{enumerate}
    \item Let $m \mid D_1N_1$ with $\omega(\gcd(D_1,m))$ and $\omega(\gcd(D_2,m))$ having opposite parity. If $w_m(D_1,N_1) \in W_1$, then $w_m(D_2,N_2) \not \in W_2$. 
    \item Let $m \mid D_1N_1$ with $\omega(\gcd(D_1,m))$ and $\omega(\gcd(D_2,m))$ having the same parity, and suppose that $g \coloneqq g(X_0^{D_1}(N_1)/W_1) = g(X_0^{D_2}(N_2)/W_2) \geq 2$. Then
    \[
    w_m(D_1,N_1) \in W_1 \text{ if and only if } w_m(D_2,N_2) \in W_2.\]
\end{enumerate}
Moreover, if $g \geq 2$ then we have the following equality of sets:
\begin{align*} &\left\{m \mid w_m(D_1,N_1) \in W_1 \textnormal{ and } w_m(D_2,N_2) \in W_2 \right\} = \\
&  \left\{m \mid \left(w_m(D_1,N_1) \in W_1 \textnormal{ or } w_m(D_2,N_2) \in W_2\right) \textnormal{ and } 2 \mid \omega(\gcd(D_1,m))-\omega(\gcd(D_2,m)) \right\} .
\end{align*}
\end{proposition}
\begin{proof}
First, note that the assumption that there exist newforms of \emph{all} sign patterns in level $D_1N_1$ implies that the Atkin--Lehner quotient $X_0^{D_i}(N_i)/W_i$ has positive genus for $i \in \{1,2\}$. The equality of sets at the end of the statement follows immediately from Parts (1) and (2) of the statement, so we proceed with their proofs.
\begin{enumerate}
    \item The assumption that $\omega(\gcd(m,D_1))$ and $\omega(\gcd(m,D_2))$ have different parity implies, via \cref{Ribet_isog}, that the actions of $w_m(D_1,N_1)$ and $w_m(D_2,N_2)$ have opposite associated signs. By our hypothesis on $D_1N_1$, we then know that there exists a nonzero $f \in S_2(D_1N_1)^\text{new}$ which is compatible with the action of $w_m(D_1,N_1)$ and not with that of $w_m(D_2,N_2)$, and moreover we can take $f$ compatible with the action of $W_1$. The abelian variety $A_f$ is therefore an isogeny factor of $\Jac\left(X_0^{D_1}(N_1)/W_1 \right)$ but is not an isogeny factor of 
    \[ \Jac\left(X_0^{D_2}(N_2)/\langle w_m(D_2,N_2)\rangle \right), \]
    hence is not an isogeny factor of $\Jac\left(X_0^{D_2}(N_2)/W_2\right)$,
    giving the claim.
    \item Assume that $g(X_0^{D_1}(N_1)/W_1) \geq 2$, and assume by way of contradiction that $w_m(D_2,N_2) \not \in W_2$. Let $ H \coloneqq \langle W_2 \cup \{w_m \} \rangle \leq W_0(D_2,N_2)$. 
    By the genus assumption, we know that there is an isogeny factor $A$ of $\Jac\left(X_0^{D_2}(N_2)/W_2\right)$ which is not an isogeny factor of $\Jac\left(X_0^{D_2}(N_2)/H\right)$. Hence, there is some sign pattern which is compatible with the action of $W_2$ but not with that of $H$. By our assumption on $D_1N_1$, we know there is a newform $f$ in level $D_1N_1$ of this sign pattern and we can take $A = A_f$. Thus, $A_f$ cannot be an isogeny factor of $\Jac\left(X_0^{D_1}(N_1)/\langle w_m(D_1,N_1) \rangle\right)$, and hence not of $\Jac\left(X_0^{D_1}(N_1)/W_1\right)$; its sign pattern is not compatible with the action of $w_m(D_2,N_2)$, and because $\omega(\gcd(D_1,N_1))$ and $\omega(\gcd(D_2,N_2))$ have the same parity it similarly is not compatible with the action of $w_m(D_1,N_1)$ by \cref{Ribet_isog}. This then contradicts the existence of the isogeny of Jacobians we assumed. The converse is obtained by swapping the roles of $W_1$ and $W_2$; given the fact that we can take the isogeny factor $A$ of choice to correspond to a newform of level $D_1N_1$, and hence to be $D_i$-new for $i \in \{1,2\}$, the argument works for both directions. \qedhere
\end{enumerate}
\end{proof}


\subsection{Atkin--Lehner quotients of genus at most one}\label{subsection: low_genus}

We focus here on Atkin--Lehner quotients $X \colonequals X_0^D(N)/W$ of genus $g \leq 1$, mainly so that we may avoid providing separate arguments for these cases in what follows. Of course, if $X$ has genus $0$ then it is geometrically isomorphic to every other genus $0$ quotient of this form, and if $D=1$ then $X$ certainly has a rational point and hence is isomorphic over $\Q$ to any other genus $0$ Atkin--Lehner quotient with a rational point. If $X$ has genus $1$ then it is geometrically isomorphic to its Jacobian, so the determination of isogenous and isomorphic quotients in genus $1$ just comes down to determining these Jacobian isomorphism classes and the existence of rational points (which, again, is automatic if $D=1$).

The determination of all such $X$ of genus $0$ and genus $1$ in the $D>1$ case was completed by Padurariu--Saia in \cite{PS25}. In the referenced work, the authors also determine all isomorphism classes of Jacobians of such $X$ of genus $1$ with $N$ squarefree, and they make progress on determining which $X$ of genus $0$ and $1$ have a rational point. We refer to their work for the $D>1$ case; we do not seek to improve on their results here. 

Rather, we devote the rest of this subsection to the $D=1$ case. The following result for the quotients $X_0(N)^* \coloneqq X_0(N)/W_0(1,N)$ can be found in work of Gonz\'{a}lez--Lario.\footnote{Note that there is a typo therein; the table in \cite[Proposition 3.2]{GL} lists $114$ twice, whereas the entry $114$ in the $\omega(N) = 2$ row should be $141$.}

\begin{proposition}\cite[Proposition 3.1, Proposition 3.2]{GL}\label{GL_prop}
\begin{enumerate}
    \item There are exactly $43$ squarefree positive integers $N>1$ such that $X_0(N)^*$ has genus $0$, namely those in the following table:
    
\begin{center}
    \begin{tabular}{|c|c|} \hline 
$\omega(N)$ & $N$ \\ \hline 
$1$ & $2, 3, 5, 7, 11, 13, 17, 19, 23, 29, 31, 41, 47, 59, 71$ \\
$2$ & $6, 10, 14, 15, 21, 22, 26, 33, 34, 35, 38, 39, 46, 51, 55, 62, 69, 87, 94, 95, 119$ \\
$3$ & $30, 42, 66, 70, 78, 105, 110$ \\ \hline
\end{tabular}
\end{center}

\item There are exactly $38$ squarefree positive integers $N$ such that $X_0(N)^*$ has genus $1$, namely those in the following table:

\begin{center}
    \begin{tabular}{|c|c|} \hline 
$\omega(N)$ & $N$ \\ \hline 
$1$ & $37, 43, 53, 61, 79, 83, 89, 101, 131$ \\
$2$ & $57, 58, 65, 74, 77, 82, 86, 91, 111, 118, 123, 141, 142, 143, 145, 155, 159$ \\
$3$ & $102, 114, 130, 138, 174, 182, 190, 195, 222, 231, 238$ \\ 
$4$ & $210$ \\ \hline
\end{tabular}
\end{center}

\end{enumerate}
\end{proposition}

While it is a relatively simple task to enumerate all genus $0$ and $1$ quotients of the form $X_0(N)/W$ from the information of the above proposition, we are not aware of a source where this is completed in the literature. Similarly, while we are able to find isomorphism classes of certain genus $1$ Atkin--Lehner quotients of modular curves (e.g., one can find many, though not all, cases of the form $X_0(N)/\langle w_m \rangle$ for some $m \parallel N$ in \cite[Tables 1 and 2]{Gonzalez_bielliptics}), we do not know of a reference that handles arbitrary Atkin--Lehner subgroups. For this reason, we provide enumerations of the genus $0$ and genus $1$ Atkin--Lehner quotients of $X_0(N)$, for $N$ squarefree, in the files \repo{LowGenusModular/genus_0_modular_quotients.m}{\texttt{genus{\_}0{\_}modular{\_}quotients.m}} and \repo{LowGenusModular/genus_1_modular_quotients.m}{\texttt{genus{\_}1{\_}modular{\_}quotients.m}}, and we determine the isomorphism classes of all genus $1$ quotients in the following proposition.

\begin{proposition}
If $N$ is a squarefree positive integer and $W \leq W_0(1,N)$ is such that $X \coloneqq X_0(N)/W$ has genus $1$ (and hence is an elliptic curve), then the isomorphism class of $X$ is given in \cref{table: genus_1_modular_quotients}. 
\end{proposition}
\begin{proof}
First of all, to be clear about our enumeration of genus $1$ (and genus $0$) quotients: if $X_0^D(N)/W$ has genus $g \leq 1$, then $X_0(N)^*$ also has genus $g^* \leq 1$ and hence $N$ is among the levels listed in \cref{GL_prop}. We arrive at the list of all genus $1$ (and that of all genus $0$) Atkin--Lehner quotients simply by computing, for each level appearing in this referenced proposition, the genus of $X_0(N)/W$ for every $W \leq W_0(D,N)$. 

Having obtained the list of all genus $1$ quotients, we determine the isomorphism class of each quotient that does not appear in \cite{Gonzalez_bielliptics} using functionality for modular curve quotients available in Magma. All relevant code can be found in the file \repo{LowGenusModular/low_genus_modular_checks.m}{\texttt{low{\_}genus{\_}checks.m}}. 
\end{proof}

\begin{remark}
The curve $X_0(51)/\langle w_{51} \rangle$, which is isomorphic to the elliptic curve with Cremona reference 17a4, is mistakenly identified as that with reference 51a1 in \cite{Gonzalez_bielliptics}. 
\end{remark}

{\rowcolors{2}{white}{gray!20}
\setlength{\tabcolsep}{2pt}
\begin{longtable}{|c|c|c||c|c|c||c|c|c|}
\caption{Isomorphism classes of quotients $X_0(N)/W$, with $N$ squarefree, of genus $1$, given by their Cremona references.}\label{table: genus_1_modular_quotients} \\ \hline

\rowcolor{gray!50}
$N$ & $W$ & Reference & $N$ & $W$ & Reference & $N$ & $W$ & Reference \\ \hline
$ 11 $ & $\langle
\textnormal{id} \rangle $
 &  11a1  & 
$ 14 $ & $\langle
\textnormal{id} \rangle $
 &  14a1  & 
$ 14 $ & $\langle
w_{ 2 } \rangle $
 &  14a2  \\ \hline
$ 15 $ & $\langle
\textnormal{id} \rangle $
 &  15a1  & 
$ 15 $ & $\langle
w_{ 3 } \rangle $
 &  15a4  & 
$ 17 $ & $\langle
\textnormal{id} \rangle $
 &  17a1  \\ \hline
$ 19 $ & $\langle
\textnormal{id} \rangle $
 &  19a1  & 
$ 21 $ & $\langle
\textnormal{id} \rangle $
 &  21a1  & 
$ 21 $ & $\langle
w_{ 7 } \rangle $
 &  21a3  \\ \hline
$ 22 $ & $\langle
w_{ 2 } \rangle $
 &  11a1  & 
$ 22 $ & $\langle
w_{ 22 } \rangle $
 &  11a3  & 
$ 26 $ & $\langle
w_{ 2 } \rangle $
 &  26a1  \\ \hline
$ 26 $ & $\langle
w_{ 13 } \rangle $
 &  26b1  & 
$ 30 $ & $\langle
w_{ 5 } \rangle $
 &  30a1  & 
$ 30 $ & $\langle
w_{ 6 } \rangle $
 &  15a3  \\ \hline
$ 30 $ & $\langle
w_{ 30 } \rangle $
 &  15a8  & 
$ 30 $ & $\langle
w_{ 2 }, 
w_{ 3 } \rangle $
 &  15a1  & 
$ 30 $ & $\langle
w_{ 2 }, 
w_{ 5 } \rangle $
 &  30a2  \\ \hline
$ 30 $ & $\langle
w_{ 3 }, 
w_{ 10 } \rangle $
 &  15a3  & 
$ 33 $ & $\langle
w_{ 33 } \rangle $
 &  11a3  & 
$ 34 $ & $\langle
w_{ 2 } \rangle $
 &  17a2  \\ \hline
$ 34 $ & $\langle
w_{ 17 } \rangle $
 &  34a1  & 
$ 34 $ & $\langle
w_{ 34 } \rangle $
 &  17a4  & 
$ 35 $ & $\langle
w_{ 5 } \rangle $
 &  35a1  \\ \hline
$ 37 $ & $\langle
w_{ 37 } \rangle $
 &  37a1  & 
$ 38 $ & $\langle
w_{ 19 } \rangle $
 &  38b1  & 
$ 38 $ & $\langle
w_{ 38 } \rangle $
 &  19a3  \\ \hline
$ 39 $ & $\langle
w_{ 3 } \rangle $
 &  39a1  & 
$ 42 $ & $\langle
w_{ 14 } \rangle $
 &  21a4  & 
$ 42 $ & $\langle
w_{ 2 }, 
w_{ 3 } \rangle $
 &  14a1  \\ \hline
$ 42 $ & $\langle
w_{ 2 }, 
w_{ 7 } \rangle $
 &  21a1  & 
$ 42 $ & $\langle
w_{ 2 }, 
w_{ 21 } \rangle $
 &  14a4  & 
$ 42 $ & $\langle
w_{ 3 }, 
w_{ 7 } \rangle $
 &  42a1  \\ \hline
$ 42 $ & $\langle
w_{ 6 }, 
w_{ 7 } \rangle $
 &  21a4  & 
$ 43 $ & $\langle
w_{ 43 } \rangle $
 &  43a1  & 
$ 51 $ & $\langle
w_{ 17 } \rangle $
 &  51a1  \\ \hline
$ 51 $ & $\langle
w_{ 51 } \rangle $
 &  17a4  & 
$ 53 $ & $\langle
w_{ 53 } \rangle $
 &  53a1  & 
$ 55 $ & $\langle
w_{ 11 } \rangle $
 &  55a1  \\ \hline
$ 55 $ & $\langle
w_{ 55 } \rangle $
 &  11a3  & 
$ 57 $ & $\langle
w_{ 3 }, 
w_{ 19 } \rangle $
 &  57a1  & 
$ 58 $ & $\langle
w_{ 2 }, 
w_{ 29 } \rangle $
 &  58a1  \\ \hline
$ 61 $ & $\langle
w_{ 61 } \rangle $
 &  61a1  & 
$ 62 $ & $\langle
w_{ 31 } \rangle $
 &  62a1  & 
$ 65 $ & $\langle
w_{ 65 } \rangle $
 &  65a1  \\ \hline
$ 65 $ & $\langle
w_{ 5 }, 
w_{ 13 } \rangle $
 &  65a2  & 
$ 66 $ & $\langle
w_{ 2 }, 
w_{ 11 } \rangle $
 &  66a1  & 
$ 66 $ & $\langle
w_{ 2 }, 
w_{ 33 } \rangle $
 &  11a3  \\ \hline
$ 66 $ & $\langle
w_{ 3 }, 
w_{ 11 } \rangle $
 &  66b1  & 
$ 69 $ & $\langle
w_{ 23 } \rangle $
 &  69a1  & 
$ 70 $ & $\langle
w_{ 2 }, 
w_{ 35 } \rangle $
 &  14a4  \\ \hline
$ 70 $ & $\langle
w_{ 5 }, 
w_{ 7 } \rangle $
 &  70a1  & 
$ 70 $ & $\langle
w_{ 5 }, 
w_{ 14 } \rangle $
 &  35a3  & 
$ 74 $ & $\langle
w_{ 2 }, 
w_{ 37 } \rangle $
 &  37a1  \\ \hline
$ 77 $ & $\langle
w_{ 7 }, 
w_{ 11 } \rangle $
 &  77a1  & 
$ 78 $ & $\langle
w_{ 2 }, 
w_{ 39 } \rangle $
 &  26a3  & 
$ 78 $ & $\langle
w_{ 3 }, 
w_{ 13 } \rangle $
 &  26b1  \\ \hline
$ 78 $ & $\langle
w_{ 3 }, 
w_{ 26 } \rangle $
 &  39a4  & 
$ 79 $ & $\langle
w_{ 79 } \rangle $
 &  79a1  & 
$ 82 $ & $\langle
w_{ 2 }, 
w_{ 41 } \rangle $
 &  82a1  \\ \hline
$ 83 $ & $\langle
w_{ 83 } \rangle $
 &  83a1  & 
$ 86 $ & $\langle
w_{ 2 }, 
w_{ 43 } \rangle $
 &  43a1  & 
$ 89 $ & $\langle
w_{ 89 } \rangle $
 &  89a1  \\ \hline
$ 91 $ & $\langle
w_{ 7 }, 
w_{ 13 } \rangle $
 &  91a1  & 
$ 94 $ & $\langle
w_{ 47 } \rangle $
 &  94a1  & 
$ 95 $ & $\langle
w_{ 95 } \rangle $
 &  19a3  \\ \hline
$ 101 $ & $\langle
w_{ 101 } \rangle $
 &  101a1  & 
$ 102 $ & $\langle
w_{ 2 }, 
w_{ 3 }, 
w_{ 17 } \rangle $
 &  102a1  & 
$ 105 $ & $\langle
w_{ 3 }, 
w_{ 35 } \rangle $
 &  15a8  \\ \hline
$ 105 $ & $\langle
w_{ 5 }, 
w_{ 7 } \rangle $
 &  21a4  & 
$ 105 $ & $\langle
w_{ 5 }, 
w_{ 21 } \rangle $
 &  35a3  & 
$ 105 $ & $\langle
w_{ 15 }, 
w_{ 21 } \rangle $
 &  105a1  \\ \hline
$ 110 $ & $\langle
w_{ 2 }, 
w_{ 55 } \rangle $
 &  11a3  & 
$ 110 $ & $\langle
w_{ 5 }, 
w_{ 11 } \rangle $
 &  110b1  & 
$ 110 $ & $\langle
w_{ 10 }, 
w_{ 11 } \rangle $
 &  55a4  \\ \hline
$ 111 $ & $\langle
w_{ 3 }, 
w_{ 37 } \rangle $
 &  37a1  & 
$ 114 $ & $\langle
w_{ 2 }, 
w_{ 3 }, 
w_{ 19 } \rangle $
 &  57a1  & 
$ 118 $ & $\langle
w_{ 2 }, 
w_{ 59 } \rangle $
 &  118a1  \\ \hline
$ 119 $ & $\langle
w_{ 119 } \rangle $
 &  17a4  & 
$ 123 $ & $\langle
w_{ 3 }, 
w_{ 41 } \rangle $
 &  123b1  & 
$ 130 $ & $\langle
w_{ 2 }, 
w_{ 5 }, 
w_{ 13 } \rangle $
 &  65a1  \\ \hline
$ 131 $ & $\langle
w_{ 131 } \rangle $
 &  131a1  & 
$ 138 $ & $\langle
w_{ 2 }, 
w_{ 3 }, 
w_{ 23 } \rangle $
 &  138a1  & 
$ 141 $ & $\langle
w_{ 3 }, 
w_{ 47 } \rangle $
 &  141d1  \\ \hline
$ 142 $ & $\langle
w_{ 2 }, 
w_{ 71 } \rangle $
 &  142b  & 
$ 143 $ & $\langle
w_{ 11 }, 
w_{ 13 } \rangle $
 &  143a1  & 
$ 145 $ & $\langle
w_{ 5 }, 
w_{ 29 } \rangle $
 &  145a1  \\ \hline
$ 155 $ & $\langle
w_{ 5 }, 
w_{ 31 } \rangle $
 &  155c1  & 
$ 159 $ & $\langle
w_{ 3 }, 
w_{ 53 } \rangle $
 &  53a1  & 
$ 174 $ & $\langle
w_{ 2 }, 
w_{ 3 }, 
w_{ 29 } \rangle $
 &  58a1  \\ \hline
$ 182 $ & $\langle
w_{ 2 }, 
w_{ 7 }, 
w_{ 13 } \rangle $
 &  91a1  & 
$ 190 $ & $\langle
w_{ 2 }, 
w_{ 5 }, 
w_{ 19 } \rangle $
 &  190b1  & 
$ 195 $ & $\langle
w_{ 3 }, 
w_{ 5 }, 
w_{ 13 } \rangle $
 &  65a1  \\ \hline
$ 210 $ & $\langle
w_{ 2 }, 
w_{ 3 }, 
w_{ 5 }, 
w_{ 7 } \rangle $
 &  210d1  & 
$ 222 $ & $\langle
w_{ 2 }, 
w_{ 3 }, 
w_{ 37 } \rangle $
 &  37a1  & 
$ 231 $ & $\langle
w_{ 3 }, 
w_{ 7 }, 
w_{ 11 } \rangle $
 &  77a1  \\ \hline
$ 238 $ & $\langle
w_{ 2 }, 
w_{ 7 }, 
w_{ 17 } \rangle $
 &  238b1  & & & & & & \\ \hline
 
\end{longtable}}


\subsection{Sign patterns admitting newforms}\label{section: PPV}

It has long been known that \(S_2(p)\) is trivial -- equivalently, that \(g(X_0(p))=0\) -- only for primes \(p=2,3,5,7,13\), and that \(\dim S_2(p)=g(X_0(p))=1\) only for primes \(p=11, 17, 19\); both follow from the genus formula. Since we study Atkin--Lehner quotients, we will also be interested in when the subspaces of cusp forms with specific sign patterns are trivial. For squarefree \(N\), the newspace decomposes as
\[
  S_2(N)^{\mathrm{new}}
  = \bigoplus_{\varepsilon \in \{\pm 1\}^{\omega(N)}} S_2(N)^{\mathrm{new},\varepsilon},
\]
where the sum goes over all sign patterns of involutions \(w_p\) for primes \(p\mid N\). 

Only finitely many pairs \((N,\varepsilon)\) with \(N\) squarefree have
\(\dim S_2(N)^{\mathrm{new},\varepsilon} = 0\). The determination of all such pairs is the content of \cite[Theorem 1.1]{PPV26}, and one can find these pairs listed in \cite[Tables 2-5]{PPV26}. In particular, it follows from their work that for all squarefree levels \(N > 4290\), there are no trivial sign pattern subspaces, i.e., that \( \dim S_2(N)^{\mathrm{new},\varepsilon} > 0\) for all sign patterns \(\varepsilon\). 

Correspondingly, there are only finitely many squarefree $N$ such that \(X_0(N)^*\) admits a prime of good reduction $p \mid N$, and more generally there are only finitely many squarefree pairs  \((D,N)\) such that \(X_0^D(N)^*\) admits a prime of good reduction $p \mid DN$. All such levels are listed in \cite[Tables 6-7]{PPV26}. For brevity, in what follows, whenever we refer to the work of Padurariu--Park--Voight \cite{PPV26}, we are specifically referring to these tables.


\subsection{Quotients of the same Shimura curve (\texorpdfstring{\ensuremath{D_1 = D_2}}{D1 = D2} and \texorpdfstring{\ensuremath{N_1 = N_2}}{N1 = N2})}\label{subsection: isogenies_same_curve}

For the entirety of this subsection, we let $D$ be the discriminant of an indefinite quaternion algebra over $\Q$, we let $N$ be a squarefree positive integer that is coprime to $D$, and we let $W_1, W_2 \leq W_0(D,N)$ be distinct Atkin--Lehner subgroups. We seek to determine the existence of isogenies of the form
\[ \Jac\left(X_0^D(N)/W_1\right) \sim \Jac\left(X_0^D(N)/W_2\right). \]
The following proposition gives a rather strong restriction on this phenomenon, as one may expect. 

\begin{proposition}\label{prop: AL_restrictions_same_curve}
Assume that there exists a newform of level $DN$ in every sign pattern. If $W_1 \ne W_2$ and
\[ \Jac\left(X_0^D(N)/W_1\right) \sim \Jac\left(X_0^D(N)/W_2\right), \]
then 
\[ g \coloneqq g\left(X_0^D(N)/W_1\right) = g\left(X_0^D(N)/W_2\right) \leq 1. \]
\end{proposition}
\begin{proof}
Taking an involution $w_m \in W_1 \setminus W_2$ (swapping $W_1$ and $W_2$ as needed), we consider the subgroup $H \coloneqq \langle W_2 \cup \{w_m\}\rangle$ and the natural quotient map
\[ X_0^D(N)/W_2 \to X_0^D(N)/H. \]
If $g \geq 2$, then it follows that 
\[ g\left(X_0^D(N)/H\right) < g. \] 
Thus, there exists a sign pattern $\overline{\epsilon}$ in level $DN$ that is compatible with the action of $W_2$ but not with that of $H$, and by the assumption on $DN$ we can take a newform $f$ in $S_2(DN)^{\overline{\epsilon}}$. The modular abelian variety $A_f$ then represents an isogeny factor of $\Jac\left(X_0^D(N)/W_2\right)$ which is not an isogeny factor of $\Jac\left(X_0^D(N)/\langle w_m \rangle\right)$ and hence is not an isogeny factor of $\Jac\left(X_0^D(N)/W_1\right)$. This contradicts the assumed isogeny of Jacobians. 
\end{proof}

\begin{remark}
Note that for general $N$ (not necessarily) squarefree, we have that
\[\Jac(X_0^D(N)/W_1) \sim_{\Q}\Jac(X_0^D(N)/W_2) \implies W_1 = W_2,\] 
but when $N$ is squarefree the same implication holds for isogenies over $\overline{\Q}$.
\end{remark}

By \cref{prop: AL_restrictions_same_curve}, we can begin by restricting our attention to those levels $M$ appearing in the tables of Padurariu--Park--Voight \cite{PPV26}. For each such $M$ and for each divisor $D \mid M$ with $\omega(D)$ even, we let $N = M/D$, range over all Atkin--Lehner subgroups $W \leq W_0(D,N)$, and record the genus of $X_0^D(N)/W$. 

For curves $X_0^D(N)$ admitting more than one quotient of a fixed genus, we
\begin{itemize}
    \item use point counts over finite fields $\F_{p^r}$, via the algorithm implemented in \cite{MPSS25}, for primes $p \leq 100$ and powers $1 \leq r \leq 3$ , to attempt to exclude the possibility of isogeny. Then, 
    \item we explicitly compute isogeny decompositions of Jacobians of remaining candidates using \cref{Ribet_isog} in order to determine which are truly isogenous. 
\end{itemize}

We provide the complete, final list of isogenous Jacobians from this subsection in \cref{table: isogenies_same_curve}. Code for all computations described in this subsection can be found in the files
\begin{align*} 
&\repo{Isogenies/SameCurve/same_curve_isogeny_genus_checks.m}{\texttt{same{\_}curve{\_}isogeny{\_}genus{\_}checks.m}} \quad \text{and}  \\ &\repo{Isogenies/SameCurve/same_curve_isogenies_narrowing.m}{\texttt{same{\_}curve{\_}isogenies{\_}narrowing.m}}.
\end{align*}


\subsection{The case of \texorpdfstring{\ensuremath{D_1N_1 = D_2N_2}}{D1N1 = D2N2} with \texorpdfstring{\ensuremath{D_1 \neq D_2}}{D1 not equal to D2}}\label{subsection: isogenies_same_level}

For the remainder of this subsection, we fix the following notation:
\begin{itemize}
    \item $D_1 \neq D_2$ are discriminants of indefinite quaternion algebras over $\mathbb{Q}$, so without loss of generality in our arguments we may assume $D_1 \nmid D_2$,
    \item $N_1$ and $N_2$ are squarefree positive integers with $D_1N_1 = D_2N_2$, 
    \item $W_i \leq W_0(D_i,N_i)$ is an Atkin--Lehner subgroup for each $i=1,2$, and 
    \item we fix a prime factorization
    \[ D_1N_1 = p_1 \cdots p_{\omega}, \]
    ordered such that $p_1 \mid D_1$ and $p_1 \nmid D_2$.
\end{itemize}

We begin with a few results that reduce our consideration of existence of isogenies of the form 
\[ \Jac\left(X_0^{D_1}(N_1)/W_1\right) \sim \Jac\left(X_0^{D_2}(N_2)/W_2\right) \]
to finitely many levels $D_1N_1$.

\begin{proposition}\label{prop: dim_zero_prop_shimura}
Let $\overline{\epsilon} = [\epsilon_1,\ldots, \epsilon_{\omega}]$ be a sign pattern that is compatible with the actions of both $W_1$ and $W_2$. Let $\widehat{\epsilon} = [\epsilon_2,\ldots,\epsilon_{\omega}]$. If there exists an isogeny 
\[ \Jac\left(X_0^{D_1}(N_1)/W_1\right) \sim  \Jac\left(X_0^{D_2}(N_2)/W_2\right), \]
then
\[ \dim\left(S_2(D_1N_1/p_1)^{D_2-\textnormal{new},\widehat{\epsilon}}\right) = \dim\left(S_2(p_2\cdots p_{\omega})^{D_2-\textnormal{new},\widehat{\epsilon}}\right) = 0 .\]
\end{proposition}
\begin{proof}
By \cref{Ribet_isog}, we have that $\Jac\left(X_0^{D_i}(N_i)\right)$ is isogenous to $\Jac\left(X_0(D_1N_1)\right)^{D_i-\textnormal{new}}$ for $i \in \{1,2\}$. By our choice of $\overline{\epsilon}$, for every form $f$ in $S_2(D_1 N_1)^{D_i-\text{new},\overline{\epsilon}}$ the associated modular abelian variety $A_f$ represents an isogeny factor of the Jacobian of $X_0^{D_i}(N_i)/W_i$. Our hypothesis that there exists an isogeny between these two Jacobians then implies that 
\[ S_2(D_1N_1)^{D_1-\text{new},\overline{\epsilon}} = S_2(D_1N_1)^{D_2-\text{new},\overline{\epsilon}}. \]
Considering the embedding
\[ S_2(D_1N_1/p_1)^{D_2-\textnormal{new},\widehat{\epsilon}} \hookrightarrow S_2(D_1 N_1)^{\overline{\epsilon}}, \]
we note that its image is generated by elements of the difference 
\[ S_2(D_1N_1)^{D_2-\text{new},\overline{\epsilon}} \setminus S_2(D_1N_1)^{D_1-\text{new},\overline{\epsilon}}, \]
which from the above is empty. This gives the claim. 
\end{proof}

\begin{corollary}\label{cor: level_restrictions}
Suppose that there exists an isogeny 
\[ \Jac\left(X_0^{D_1}(N_1)/W_1\right) \sim \Jac\left( X_0^{D_2}(N_2)/W_2\right) \]
with $N_1$ and $N_2$ squarefree. Further assume that $g \coloneqq g\left(X_0^{D_1}(N_1)/W\right) \geq 1$. The following statements then hold. 
\begin{enumerate}
    \item Every prime divisor $p \mid D_1N_1$ is in the set
\[ \{2,3,5,7,11,13,17,19,23,29,31,37,41,43,47,53,59,61,71,79,83,89,101,131\}. \]
    \item If $\omega(D_1N_1) \geq 3$, then $D_1N_1/p_1$ must appear in the tables of Padurariu--Park--Voight. In particular, we have $\omega(D_1N_1) \leq 6$.
\end{enumerate}
\end{corollary}
\begin{proof}
For both parts, we let $\overline{\epsilon} = [\epsilon_1,\ldots,\epsilon_\omega]$ be a sign pattern that is compatible with the actions of $W_1$ and $W_2$ as in \cref{prop: dim_zero_prop_shimura}. Such a sign pattern exists by our assumption that $g$ is positive.
\begin{enumerate}
\item First, we claim that for any prime $p_i \mid D_1N_1$ we must have
\[ \dim S_2(p_i)^{[\epsilon_i]} = 0. \]
This is because, by our assumption that $D_1 \neq D_2$, there exists at least one prime $p_j \neq p_i$ with $p_j$ dividing exactly one of $D_1$ or $D_2$. If, say, $p_j \mid D_2$ and $p_j \nmid D_1$, then we have embeddings 
\[ S_2(p_i)^{[\epsilon_i]} \hookrightarrow S_2(D_1N_1/p_j)^{[\epsilon_1,\ldots,\epsilon_{j-1},\epsilon_{j+1},\ldots,\epsilon_{r}]} \hookrightarrow  S_2(D_1N_1)^{D_2-\text{new},\overline{\epsilon}} \setminus S_2(D_1N_1)^{D_1-\text{new},\overline{\epsilon}}, \]
and by \cref{prop: dim_zero_prop_shimura} the space in the middle must be zero-dimensional. (If $p_j \mid D_1$ and $p_j \nmid D_2$, one need only swap ``${D_1-\text{new}}$'' and ``${D_2-\text{new}}$'' here, achieving the same conclusion.)

We now prove the statement of Part (1). If $\epsilon_i = 1$, then we have
\[ 0 = \dim S_2(p_i)^{[+]} = \dim \Jac\left(X_0(p_i)^* \right) = g(X_0(p_i)^*). \] 
If $\epsilon_i = -1$, then the condition
\[ 0 = \dim S_2(p_i)^{[-]} \]
implies that 
\[ g(X_0(p_i)) = \dim S_2(p_i) = \dim S_2(p_i)^{[+]} = g(X_0(p_i)^*), \]
which forces $g(X_0(p_i)) \leq 1$. In either case, the claim follows from \cref{GL_prop}.

\item Letting $p_1$ and $\widehat{\epsilon}$ be as in \cref{prop: dim_zero_prop_shimura}, i.e., $p_1 \mid D_1,$ $p_1 \nmid D_2$, and $\widehat{\epsilon} = [\epsilon_2, \dotsc, \epsilon_\omega]$, we must have by this referenced proposition that
\[  \dim\left(S_2(DN/p_1)^{D_2-\textnormal{new},\widehat{\epsilon}}\right) = 0. \]
In particular, we have
\[ \dim\left(S_2(DN/p_1)^{\textnormal{new}, \widehat{\epsilon}}\right) = 0. \]
By \cite{PPV26}, we then know that $DN/p_1$ is among the finitely many levels lacking newforms for $\overline{\epsilon}$, and in particular $\omega(DN/p_1) \leq 5$ such that $\omega(DN) \leq 6$.  
\end{enumerate}
\end{proof}

The following proposition is proven using the same argument seen in the proof of \cref{prop: dim_zero_prop_shimura}, but is generalized to be useful in excluding the possibility of isogeny more widely.

\begin{proposition}\label{prop: general_dim_zero_prop_shimura}
    Let $D$ be a divisor of $D_1N_1$ such that $D_1 \mid D$ and $D_2 \nmid D$. Let $\overline{\epsilon}$ be a sign pattern in level $D_1N_1$ which is compatible with the action of $W_1$, and assume that
     \[ \dim\left(S_2(D)^{\textnormal{new},\overline{\epsilon}}\right) > 0. \] 
     Then $\Jac(X_0^{D_1}(N_1)/W_1)$ and $\Jac(X_0^{D_2}(N_2)/W_2)$ are not isogenous.
\end{proposition}
\begin{proof}
    There exists an embedding
    \[S_2(D)^{\textnormal{new},\overline{\epsilon}} \hookrightarrow S_2(D_1N_1)^{D_1-\textnormal{new},\overline{\epsilon}} \setminus S_2(D_1N_1)^{D_2-\text{new}}. \]
    Taking a newform $f$ in the image of this embedding, we have that $A_f$ is an isogeny factor of $\Jac(X_0^{D_1}(N_1)/W_1)$ but not of $\Jac(X_0^{D_2}(N_2)/W_2)$. 
\end{proof}

By \cref{cor: level_restrictions}, we may begin by restricting our attention to levels of the form $M \cdot p$, where $p$ is a prime listed in Part (1) of the referenced corollary and $M$ is a level lacking newforms for at least one sign pattern. For all such levels $M \cdot p$, we compute the genera of all quotients of the form $X_0^D(N)/W$ where $DN = M \cdot p$ and $W \leq W_0(D,N)$ and record all instances of genus matches in this level. 

Of course, equality of genera of $X$ and $Y$ is a coarse check on the existence of an isogeny $\Jac(X) \sim \Jac(Y)$. For all genus matches in each specified level, we next run through the following sequence of additional checks to exclude isogeny candidates:
\begin{itemize}
    \item we use \cref{prop: AL_restrictions} with the coarse information of which levels have at least one sign pattern lacking newforms,
    \item we use \cref{prop: general_dim_zero_prop_shimura} with the finer information of which specific sign patterns in a given level lack newforms, and
    \item we use point counts of Atkin--Lehner quotients over finite fields $\F_{p^r}$, via the algorithm implemented in \cite{MPSS25}, for primes $p \leq 100$ and powers $1 \leq r \leq 3$. 
\end{itemize}
Lastly, now that we have reduced to a reasonable number of candidates of reasonable level, we explicitly compute isogeny decompositions of Jacobians of remaining candidates using \cref{Ribet_isog} and determine those which are truly isogenous. We provide the complete, final list of isogenous Jacobians from this subsection in \cref{table: isogenies_same_product}. Code for all computations described in this subsection can be found in the files
\begin{align*} 
&\repo{Isogenies/SameProductLevel/same_product_isogeny_genus_checks.m}{\texttt{same{\_}product{\_}isogeny{\_}genus{\_}checks.m}} \quad \text{and}  \\ &\repo{Isogenies/SameProductLevel/same_product_isogenies_narrowing.m}{\texttt{same{\_}product{\_}isogenies{\_}narrowing.m}}.
\end{align*}


\subsection{The case of \texorpdfstring{\ensuremath{D_1N_1 \neq D_2N_2}}{D1N1 not equal to D2N2}}\label{subsection: isogenies_different_level}

For the remainder of this subsection, we take the following notation:
\begin{itemize}
    \item $D_1$ and $D_2$ are discriminants of indefinite quaternion algebras over $\mathbb{Q}$,
    \item $N_1$ and $N_2$ are squarefree positive integers with $\gcd(D_i, N_i)=1$ (for $i=1,2$) such that $D_1N_1 \neq D_2N_2$, so without loss of generality we may assume that $D_1N_1 \nmid D_2N_2$ in our arguments, and
    \item $W_i \leq W_0(D_i,N_i)$ is an Atkin--Lehner subgroup for each $i=1,2$.
\end{itemize}

We begin with a few results that will reduce our consideration of existence of isogenies of the form 
\[ \Jac\left(X_0^{D_1}(N_1)/W_1\right) \sim \Jac\left(X_0^{D_2}(N_2)/W_2\right) \]
to finitely many levels $D_1N_1$.

\begin{proposition}\label{level_change_prop}
Let $\overline{\epsilon}$ be a sign pattern in level $D_1N_1$ which is compatible with the action of $W_1$. If there is an isogeny 
\[ \Jac\left(X_0^{D_1}(N_1)/W_1\right) \sim \Jac\left(X_0^{D_2}(N_2)/W_2\right) \]
then 
\[ \dim\left(S_2(D_1N_1)^{\textnormal{new}, \overline{\epsilon}}\right) = 0. \]
\end{proposition}
\begin{proof}
Assume by way of contradiction that $\dim(S_2(D_1N_1)^{\overline{\epsilon},\textnormal{new}})>0$. Choose a newform $f$ of level $D_1N_1$ with sign pattern $\overline{\epsilon}$, so that $A_f$ represents an isogeny factor of $\Jac(X_0^{D_1}(N_1)/W_1)$. Every simple isogeny factor of $\Jac(X_0^{D_2}(N_2))$ is isogenous to the modular abelian variety associated to a newform $g$ of level dividing $D_2N_2$, and so $A_f$ does not represent an isogeny factor of this Jacobian.
\end{proof}

\begin{corollary}\label{cor:distinct_products}
Assume the same setup as in \cref{level_change_prop}, with $X_0^{D_1}(N_1)/W_1$ of positive genus. Further suppose that $N_1$ is squarefree, and let $p\parallel D_1N_1$ be a prime with $p\nmid D_2N_2$. Then $\Jac(X_0^{D_1}(N_1)/W_1)$ has good reduction at $p$, and hence so does $\Jac(X_0^{D_1}(N_1)^*)$. Furthermore, $X_0^{D_1}(N_1)^*$ has good reduction at $p$ as well.
\end{corollary}
\begin{proof}
Every simple isogeny factor of $\Jac(X_0^{D_2}(N_2)/W_2)$ is $A_g$ for a newform $g$ of level dividing $D_2N_2$. Since $p\nmid D_2N_2$, $\Jac(X_0^{D_2}(N_2)/W_2)$ has good reduction at $p$, and hence $\Jac(X_0^{D_1}(N_1)/W_1)$ and $\Jac(X_0^{D_1}(N_1)^*)$ both have good reduction at $p$. In particular, every newform contributing to $\Jac(X_0^{D_1}(N_1)/W_1)$ has level coprime to $p$. As the genus is positive, such a newform exists, and its level is divisible by $D_1$. Hence $p\mid N_1$. Good reduction of $X_0^{D_1}(N_1)^*$ at $p$ follows from \cite[Proposition 3.3]{PPV26}.
\end{proof}

By \cref{cor:distinct_products}, we may now restrict our attention to levels dividing those appearing in the Padurariu--Park--Voight tables recording gains of primes of good reduction for star quotients, or dividing $DN$ for pairs \( (D, N)\) such that \( g(X_0^D(N)^*)=0 \). Such pairs \( (D, N)\) were determined in \cite{PS25} for $D>1$, while those with $D=1$ are listed in \cref{GL_prop}. For each such level \( L =DN \), we consider every factorization with \( \omega(D) \) even and compute the genera of all quotients \( X_0^D(N)/W\), with $W \leq W_0(D,N)$, of positive genus.

Before recording a pair of curves of the same genus as a candidate for isogeny, we check the restrictions coming from
\cref{level_change_prop} and \cref{cor:distinct_products}. In particular, when
one total level does not divide the other, we require that the corresponding
Atkin--Lehner quotient have no compatible new contribution at the missing
primes. This can be easily checked computationally for a given level, or one can refer to the tables of \cite{PPV26}. These filters lead us to a finite list of genus-compatible candidates for isogeny.

For these remaining candidates, we next use point counts of Atkin--Lehner quotients
over finite fields \(\F_{p^r}\), with similar parameters as previously (\(p < 100, 1 \leqslant r \leqslant 5\)). We are left with \(3599\) candidates
\[
\Jac\left(X_0^{D_1}(N_1)/W_1\right) \sim
\Jac\left(X_0^{D_2}(N_2)/W_2\right)
\]
with \(D_1N_1 \neq D_2N_2\) and \(D_1D_2>1\), i.e., with at least one side a
quotient of a Shimura curve. The computations of the last three paragraphs are carried out in \repoDistinctProducts{distinct_products_isogeny_filters.m}.

Further on, we prove all $3599$ candidates (with $D_1D_2>1$) indeed have isogenous Jacobians. We treat quotients of classical modular curves ($D_1=D_2=1$) at the end of this subsection. 

In principle, one could now compute the full isogeny decompositions of the Jacobians of all remaining candidates using Ribet's isogeny (\cref{Ribet_isog}) and modular symbols. In practice, this computation is too slow for the number of remaining pairs. We therefore use two simpler criteria which prove the required isogenies without computing the full Jacobian decompositions. 
Both criteria are based on the same idea:
after the compatibility conditions on Atkin--Lehner signs are fixed at the lower level, there is at most one possible extension of these signs to the higher level, and the genus equality then forces the compatible old part to account for the whole Jacobian.

\begin{definition}
Let \( p \nmid D_1N_1 \) be prime, \( D_2N_2 = pD_1N_1\) squarefree, and let \( W_1 \leq W_0(D_1, N_1)\), \(W_2 \leq W_0(D_2, N_2)\). We say that \emph{\(W_2\) extends \(W_1\) at p} if restriction to the primes dividing \(D_1N_1\) is a bijection from the sign patterns in level \( D_2N_2\) compatible with the action of \(W_2\) onto those in level \(D_1N_1\) compatible with the action of \(W_1\).
\end{definition}

\begin{remark}
This definition means that the compatibility conditions imposed by \(W_2\) at the primes
dividing \( D_1N_1 \) are exactly those imposed by \(W_1\), and that, once these signs
are fixed, compatibility with the action of \( W_2 \) also determines the sign at $p$. It implies the following:
\begin{enumerate}
    \item a sign pattern in level \(D_2N_2\) compatible with the action of \(W_2\) has restriction to the primes dividing \(D_1N_1\) compatible with the action of \(W_1\), and
    \item every sign pattern compatible with the action of \(W_1\) extends to exactly one sign pattern compatible with the action of \(W_2\).
\end{enumerate}
\end{remark}

\begin{lemma}\label{lem:p-old-criterion}
Let \(p \nmid DN\) be prime, let  \(W \leq W_0(D,N)\) and \( W' \leq W_0(D,Np)\), and put
\[
X = X_0^D(N)/W, \quad 
Y = X_0^D(Np)/W'.
\]
If \(W'\) extends \(W\) at \(p\), then \( \Jac(Y)^{p-\text{old}} \sim \Jac(X)\).

In particular, if \(g(X)=g(Y)\), then \( \Jac(Y)\sim \Jac(X) \).
\end{lemma}

\begin{proof}
Since $p\nmid DN$, the two degeneracy maps from level $N$ to level $Np$
give the $p$-old subvariety
\[
\Jac(X_0^D(Np))^{p\text{-old}}.
\]
On the modular side, this is the direct sum of two copies of the space of
forms of level dividing $DN$. The Atkin--Lehner involutions away from $p$
act on these two copies through the same signs as at level $DN$, while the
involution $w_p$ separates the two-dimensional old space attached to each
form into its two $w_p$-eigenspaces.

By assumption, the conditions imposed by $W'$ away from $p$ are exactly
those imposed by $W$, and the remaining condition determines the sign at
$p$. Hence, for every $W$-compatible form at level $DN$, the corresponding
$p$-old space at level $DNp$ contributes exactly one $W'$-compatible line.
Thus the $W'$-compatible $p$-old factor at level $DNp$ is isogenous
to the $W$-compatible factor at level $DN$. Therefore
\[
\Jac(Y)^{p\text{-old}}\sim\Jac(X).
\]

If $g(Y)=g(X)$, then \( \dim \Jac(Y)^{p\text{-old}} = g(X) = g(Y) = \dim \Jac(Y)\), so the compatible $p$-new factor
has dimension zero and \( \Jac(Y)\sim\Jac(X) \).
\end{proof}

 \begin{lemma}\label{lem:change_in_D_criterion}
 Let \(p \nmid D_1N_1\) be prime and \(D_2N_2 = pD_1N_1\), and let \(W \leq W_0(D_1,N_1)\) and \( W' \leq W_0(D_2,N_2)\). If \(W'\) extends \(W\) at \(p\), there is a unique embedding
\[ S_2(D_1N_1)^{W_1,D_1-\text{new}} \hookrightarrow S_2(D_2N_2)^{W'}. \]
Furthermore, \(\Jac(X_0^{D_1}(N_1)/W) \sim \Jac(X_0^{D_2}(N_2)/W')\) if and only if the image of this embedding is $S_2(D_2N_2)^{W',D_2-\text{new}}$.
\end{lemma}

\begin{proof}
The proof is analogous to the proof of \cref{lem:p-old-criterion}.
\end{proof}

\begin{remark}
In the implementation of both lemmas (in \repoDistinctProducts{provers1.m}), only the newspaces at the lower product $D_1N_1$ are computed (and the uniqueness of the extensions of sign patterns is checked). Since candidates have already been filtered so that $g(X)=g(Y)$, once the compatible newsubspace has dimension
$g(X)$, it fills the whole compatible part for $Y$, and the compatible
$p$-new subspace is necessarily zero.

Before applying the lemmas, \repoDistinctProducts{provers1.m} decides genus $1$ pairs in which both curves are quotients $X_0^{D_i}(N_i)/W_i$ with $D_i>1$ and $W_i$ nontrivial, by comparing the isogeny classes of their Jacobians determined in \cite{PS25}.
\end{remark}

Furthermore, performing a transitive closure on the so far proven isogenies reduces the number of remaining pairs to \( 661\), on which we finally run the full Jacobian decomposition (in \repoDistinctProducts{distinct_levels_direct_final_candidates.m}). All $661$ remaining pairs turn out to have isogenous Jacobians, so all $3599$ candidates do.

For the modular case $D_1=D_2=1$, we restrict to levels \( N\) dividing those in the Padurariu--Park--Voight table (\cite[Table 6]{PPV26}), together with those \( N \) such that \( g(X_0(N)^*)=0 \). For each such level we compute the genera of all quotients \( X_0(N)/W \) of positive genus, and consider every pair of distinct levels. We apply \cref{level_change_prop}, compare the genera and apply \cref{cor:distinct_products}. Finally, we compute the full Jacobian decompositions of the remaining \( 739 \) candidate pairs (at the end of \repoDistinctProducts{distinct_levels_for_isogenies_of_modular_quotients.m}); exactly \( 150 \) of them have isogenous Jacobians.
This proves the following.
\begin{corollary} The following holds.
\begin{enumerate}
    \item There are \( 643 \) pairs of quotients \( X_0^{D_1}(N_1)/W_1 \) and \( X_0^{D_2}(N_2)/W_2 \) of genus greater than \( 1 \) with isogenous Jacobians, where \( D_1N_1 \) and \( D_2N_2 \) are distinct and squarefree, and \( D_1D_2 > 1 \). All these pairs are listed in \cref{tab:isog_shimura_different_levels}.

\item \raggedright There are \( 61\) pairs of quotients \( X_0(N_1)/W_1\) and \(X_0(N_2)/W_2\) of genus greater than \( 1\) with isogenous Jacobians for squarefree \( N_1 \neq N_2 \), also listed in \cref{tab:isog_shimura_different_levels}.
\end{enumerate}
\end{corollary}


\section{Reduction of Shimura curve quotients at primes dividing the level}\label{section: reductions_at_N}

\subsection{Reduction type computations}
In \cite{PS25}, the authors implement an algorithm to compute information on the reduction types of Atkin--Lehner quotients $X_0^D(N)/W$, with $N$ squarefree, at primes $p \mid D$. This is done using the theory of {\v{C}}erednik--Drinfeld uniformizations of the Shimura curves $X_0^D(N)$. 

For primes $p \parallel N$, the reduction types of $X_0^D(N)$ are handled analogously to the theory of Deligne--Rapoport that applies to the $D=1$ case of the modular curves $X_0(N)$. Indeed, Deligne--Rapoport remark in their original work that their theory applies in the $D>1$ case as well \cite[p. 155]{DeligneRapoport}. One finds the relevant statements, which we collect below, worked out independently in work of Helm \cite{Helm07} and of Clark \cite{Clark09}. Their work follows important work of Buzzard on integral models of Shimura curves, though his work focuses on moduli problems with additional rigidifying information compared to the coarse moduli spaces that we treat in this work.

For the general background theory on reduction of algebraic curves, we refer the reader to \cite[Chapters 9 and 10]{Liu02}; more specifically, for the definition of thickness, see \cite[Corollary~10.3.22, Definition~10.3.23]{Liu02} (or \cite[Definition 2.1]{ADK26}), for semi-stable and stable models, and split nodes, see \cite[§10.3.1]{Liu02}, and for regular models and minimal desingularization of semi-stable curves, see ~\cite[Corollary 10.3.25]{Liu02}.

We fix the setup for the entire subsection.
\begin{setup}\label{setup:al_quotient} Let $p$ be a prime, let $M$ be an integer (not necessarily squarefree) coprime to $p$, set $N \colonequals p M$, let $D \geq 1$ a discriminant of an indefinite quaternion algebra over $\Q$ with $\gcd(D,N)=1$, and let $W \leq W_0(D, N)$ be a subgroup of the Atkin--Lehner operators on $X_0^D(N)$. Put $W' \colonequals \{w_m \in W : p \nmid m\}$, viewed as a subgroup of $W_0(D,M)$.
\end{setup}

We write points of $X_0^D(N)$ as $(A,G)$, suppressing the QM-structure $\iota$ from the notation; here $G\subseteq A[N]$ is $\mathcal O_B$-stable of order $N^2$, and at $p$ we allow, as in \cite{Buzzard97}, any $\mathcal O_B$-stable finite flat subgroup scheme of order $p^2$. The two degeneracy maps $X_0^D(N) \to X_0^D(M)$, where $N=pM$ with prime $p\nmid M$, are given by the following maps
\[ \pi_1 \colon (A,G) \mapsto (A, G[M]), \qquad \pi_2 \colon (A, G) \mapsto \left(A/G[p], (G[M]+G[p])/G[p]\right), \]
where $A/G[p]$ carries the QM-structure induced by $\iota$. Over $\F_p$ there are two maps
\begin{align*}
  s_1 \colon X_0^D(M)_{\F_p} &\to X_0^D(N)_{\F_p}, &
  (A, G) &\mapsto (A,\, G + \ker \Frob_p),\\
  s_2 \colon X_0^D(M)_{\F_p} &\to X_0^D(N)_{\F_p}, &
  (A, G) &\mapsto \left(A^{(p)},\, G^{(p)} + \ker \Ver_p\right),
\end{align*}
where $\Frob_p \colon A \to A^{(p)}$ is the relative Frobenius and $\operatorname{Ver}_p \colon A^{(p)} \to A$ is the Verschiebung, and $A^{(p)}$ carries the QM-structure $\iota^{(p)}$. Note that $\pi_1 \circ s_1$ and $\pi_2 \circ s_2$ are the identity. Over a field of characteristic $p$, we call a point on $X_0^D(N)$ induced by $(A,G)$ \textbf{supersingular} if $A$ is isomorphic to the square of a supersingular elliptic curve, and we call the point \textbf{ordinary} otherwise.

The following is the analogue for $D>1$ of the Deligne--Rapoport description used in the case $D=1$ (see \cite{ADK26}); the integral model is the one constructed by Buzzard \cite{Buzzard97}.

\begin{proposition}\label{prop: semistable_XD}
With notation as above (in \Cref{setup:al_quotient}), the following hold.
\begin{enumerate}
    \item The curve $X_0^D(N)$ has a model over $\Z[\frac{1}{DM}]$ which is semi-stable, and smooth outside of the supersingular points in characteristic $p$.
    \item $X_0^D(N)_{\F_p}$ consists of two irreducible components. The map $s_1$ gives an isomorphism of $X_0^D(M)_{\F_p}$ with the first component and $s_2$ one with the second component. These two components meet transversally along the supersingular points, where a supersingular point $x \in X_0^D(M)(\Fpbar)$ on the first copy is identified with $\Frob_p\, x$ on the second copy. An element $w_m \in W_0(D,N)$ sends each component to itself if and only if $p \nmid m$.
    \item If $(A,G)$ represents a supersingular point $x \in X_0^D(N)(\overline{\F}_p)$, then the thickness of the node $x$ is $\frac12 \#\Aut(A,G)$.
\end{enumerate}
\end{proposition}

\begin{proof}
Items (1) and (3) follow directly from \cite[Theorem~10.3]{Helm07} and \cite[Theorem~13]{Clark09}, as well as everything in (2) except the action of $w_m$ (i.e., the last sentence). It suffices to consider $w_q$ for primes $q \mid DN$. From the moduli description, $w_p$ swaps $\ker \Frob_p$ and $\ker \Ver_p$, hence the two components. For $q\neq p$, $w_q$ is induced by an isogeny of degree $q^2$: the quotient by $G[q]$ if $q\mid M$, and by $A[J_q]$ if $q\mid D$, where $J_q$ is the two-sided ideal of $\mathcal{O}$ of reduced norm $q$. Such an isogeny commutes with $\Frob_p$ and $\Ver_p$, so $w_q$ fixes each component.
\end{proof}

We now use the following result of Raynaud to transition, as done in \cite{ADK26} in the $D=1$ case, to general Atkin--Lehner quotients. 

\begin{proposition} \label{proposition: Raynaud}
Let $S$ be the spectrum of a discrete valuation ring $R$ with closed point $s$. 
Let $X$ be a semi-stable quasi-projective curve over $S$, endowed with 
the action of a finite group $G$. Then the quotient scheme 
$Y = X/G$ is semi-stable. More precisely, let $x \in X_s$ be a closed 
point, $y$ its image in $Y_s$. Then we have the following properties:
\begin{enumerate}
    \item If $X$ is smooth at $x$, then $Y$ is smooth at $y$;
    \item If $x$ is an ordinary double point of $X_s$, then $y$ is a 
    smooth or ordinary double point; 
    \item Let $I'$ be the image of the inertia group $I$ at $x$ in 
    $\mathrm{Aut}_R(\mathcal{O}_{X,x})$. If $x$ is of 
    thickness $m$, and if $y$ is a double point, then $y$ is of 
    thickness $mn$, where $n$ is the order of $I'$. Furthermore if $x$ is split, then $y$ is split.
    \item If $x$ is an ordinary double point then $y$ is a smooth point if and only if there is an element $\iota$ in the inertia group $I$ at $x$ that swaps the two branches at $x$.
\end{enumerate}    
\end{proposition}
\begin{proof}
This is a slightly refined version of \cite[Proposition~10.3.48]{Liu02}, already discussed in \cite[Proposition~3.2]{ADK26}.
\end{proof}

\begin{remark}
Since we have a notion of a point on a model of a Shimura curve being ordinary (i.e., not supersingular), further on we always refer to ordinary double points as \emph{nodes} in order to avoid confusion.
\end{remark}

\begin{theorem}\label{thm: semistable_XD_quotient}
As in \Cref{setup:al_quotient}, let $W' \coloneqq \{ w_m \in W : p \nmid m \}$. Then the following hold.
\begin{enumerate}
    \item The quotient $X_0^D(N)/W$ has a model $\mathcal X$ over $\Z[\frac{1}{DM}]$ which is semi-stable, and smooth outside of the images of the supersingular points in characteristic $p$.
    \item If $W' \neq W$, then $\mathcal X_{\overline{\F}_p}$ consists of a single component whose desingularization is isomorphic to $(X_0^D(M)/W')_{\overline{\F}_p}$. If $W' = W$, it consists of two components, both isomorphic to $(X_0^D(M)/W')_{\overline{\F}_p}$.
    \item Let $(A,G)$ represent a supersingular point $x \in X_0^D(N)(\overline{\F}_p)$, and let $W_x$ denote the stabilizer of $x$ in $W$. If some $w_{pd}$ lies in $W_x$, then the image of $x$ is a smooth point. Otherwise, the image of $x$ is a node of thickness $\#W_x \cdot \frac12 \#\Aut(A,G)$.
    
    \item If $g((X_0^D(N)/W)_\Q) \geqslant 2$, then $\mathcal X$ is stable.
\end{enumerate}
\end{theorem}

\begin{proof}
The proof is the same as in the modular case ($D=1$) \cite[Theorem~3.5]{ADK26}: (1) follows from \Cref{prop: semistable_XD} and \Cref{proposition: Raynaud}. The description in (2) follows from the analogous description of $X_0^D(N)_{\F_p}$, as $W'$ is the subgroup of $W$ that maps each component to itself. The description holds for all $p$, including $p=2$ (it follows from $\pi_1\circ s_1=\id$ by comparing degrees). Item (3) follows by applying \Cref{proposition: Raynaud}, using that $W_x$ acts faithfully on $\mathcal O_{X_0^D(N),x}$: if $w_m\in W_x$ acts trivially on it, then it acts trivially on its field of fractions, which is the function field of $X_0^D(N)$, so $w_m=1$ (i.e., $I'\cong W_x$ in the above notation). For (4), if there is only one component in the special fiber, then it is of (arithmetic) genus at least $2$, and if there are two components, then they either have genus $\geqslant 1$ (and an intersection point) or they are both $\P^1$ and have sufficiently many intersections (the arithmetic genus of the special fiber is then $2\cdot 0 + k-1 \geqslant 2$, where $k$ is the number of intersections).
\end{proof}

To describe the minimal regular model from the stable model, we need the two following well-known results.

\begin{corollary}[{\protect{\cite[Corollary 10.3.25]{Liu02}}}]\label{cor:regular_from_semistable_model}
Let $R$ be a DVR with closed point $s$.    Let $X \to \Spec R$ be a semi-stable projective curve 
with smooth generic fibre $X_\eta$. Let $\pi\colon X' \to X$ be the minimal desingularisation, 
$x \in X_s$ a split ordinary double point, of thickness $i$ in $X$. Then $\pi^{-1}(x)$ is 
made up of a chain of $i-1$ projective lines over $k(s)$ that meet transversally at rational 
points. 
\end{corollary}
\begin{lemma}\label{lem:minimal_regular_from_stable_model}
    Let $\mathcal X$ be a stable model over a discrete valuation ring $R$. Let $K$ be the field of fractions of $R$ and $X := \mathcal X_K$. Let $\mathcal X'$ be the minimal desingularization of $\mathcal X$. Then $\mathcal X'$ is the minimal regular model of $X$, and in particular $\mathcal X'$ is semi-stable.
\end{lemma}

\begin{proof}
This is a variation of \cite[Theorem~10.3.34(a)]{Liu02}, already established in \cite[Lemma~2.3]{ADK26}.
\end{proof}

Note that \Cref{thm: semistable_XD_quotient} part (4) and \Cref{lem:minimal_regular_from_stable_model} imply \Cref{thm: stable}.

\begin{remark}\label{rem: comparison of Frob with wp} The actions of $w_p$ and $\Frob$ (the relative Frobenius of $X_0^D(N)_{\Fpbar}$ over $\Fp$, acting on $\Fpbar$-points) on the semi-stable model of $X_0^D(N)$ can be compared as follows:
\begin{itemize}
    \item component-wise: $w_p$ swaps the two components, while $\Frob$ maps each component to itself;
    \item (crossings)-point-wise: $w_p$ and $\Frob$ agree;
    \item oriented-edges-wise: $w_p = - \Frob$.
\end{itemize}
Since $w_p$ swaps the components, denoting the components by $1$ and $2$ and taking an oriented edge $e_x$ from $1$ to $2$ corresponding to a node $x$,  this can be illustrated as in \cref{figure: Frob_action}.
{\begin{figure}[h]
\centering
\begin{tikzpicture}[thick, >=stealth,
    pt/.style={circle, fill, inner sep=1.5pt}]
  \node[pt, label=left:{$(x,1)$}]        (a) at (0,2) {};
  \node[pt, label=right:{$(\Frob \, x, 2)$}] (b) at (4,2) {};
  \node[pt, label=left:{$(x,2)$}]        (c) at (0,0) {};
  \node[pt, label=right:{$(\Frob \, x, 1)$}] (d) at (4,0) {};

  \draw[->] (a) -- node[above] {$e_x$} (b);
  \draw[->] (c) -- node[below] {$-e_{\Frob \, x}$} (d);

  \draw[|->, shorten <=6pt, shorten >=6pt] (a) -- node[left]  {$w_p$} (c);
  \draw[|->, shorten <=6pt, shorten >=6pt] (b) -- node[right] {$w_p$} (d);
\end{tikzpicture}
\caption{Action of $\Frob$ and $w_p$ on oriented edges: $w_p(e_x) = -e_{\Frob x}$}\label{figure: Frob_action}
\end{figure}
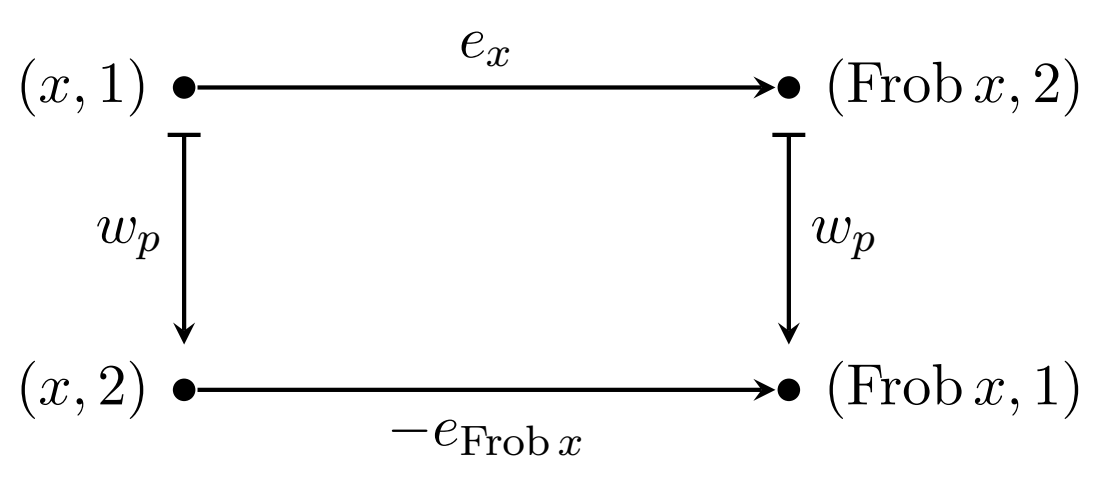}
\end{remark}

\begin{corollary}\label{cor:minimal_regular_AL_quot}
Assume $g((X_0^D(N)/W)_\Q) \geq 2$. Then the minimal regular model of $X_0^D(N)/W$ over $\Z[1/DM]$ is the minimal desingularization of $\mathcal X$, whose geometric special fiber at $p$ is obtained from the geometric special fiber $\mathcal X_{\Fpbar}$ of our semi-stable model by replacing each node $x$, the image of $(A,G)$, by a chain of $\#W_x\,\frac12\#\Aut(A,G)-1$ copies of $\P^1$. Frobenius maps the chain at $x$ to the chain at $\Frob \, x$. If $\Frob \, x = x$, it maps each $\P^1$ in the chain to itself if $x$ is split, and it reverses the chain if $x$ is non-split.
\end{corollary}
\begin{remark}
    Here $W_x$ denotes the stabilizer in $W$ of any preimage of $x$ in $X_0^D(N)$. Since $W$ is abelian, all preimages have the same stabilizer, so $W_x$ is well defined.
\end{remark}
\begin{proof}
This follows from \Cref{lem:minimal_regular_from_stable_model}, since the model $\mathcal X$ is stable (by \cref{thm: semistable_XD_quotient}(4)) and its double points have the stated thicknesses, which are replaced by chains of $\P^1$ as in \Cref{cor:regular_from_semistable_model}. As in the proof of \cite[Corollary~3.8]{ADK26}, the minimal desingularization commutes with base change to $\Z_{p^2}$, over which all nodes are split, so $\Frob$ maps the ends of the chain at $x$ to the ends of the chain at $\Frob \, x$. These ends of the chain are exactly preimages of $x$ under the normalization. If $x$ is split, this means that the ends of this chain are fixed (as they are defined over $\Fp$), which forces every $\P^1$ in the chain to be sent to itself; if $x$ is non-split, this means that the ends are swapped (as they are not defined over $\Fp$), which implies that the chain is reversed.
\end{proof}

\begin{remark}\label{rem: genus1}
Our description also holds if $g((X_0^D(N)/W)_\Q) = 1$ since the minimal desingularization of $\mathcal X$ has no $(-1)$-curves. In particular, the special fiber at $p$ of the minimal regular model is reduced.
\end{remark}

\begin{corollary}\label{cor: existence_of_Qp_points}
In the setting of \Cref{thm: semistable_XD_quotient}, $X_0^D(N)/W$ has a $\Q_p$-point if and only if at least one of the following holds:
\begin{itemize}
    \item there is a point in $(X_0^D(M)/W')(\F_p)$ that is not supersingular, or that is supersingular and has image in $\mathcal X_{\F_p}$ 
    being a smooth point or a node of thickness at least $2$;
    \item $\mathcal X_{\F_p}$ has a non-split $\F_p$-rational node of even thickness.
\end{itemize}
\end{corollary}

\begin{proof}
Let $\widetilde{\mathcal X}$ be the minimal desingularization of $\mathcal X$. It is regular and proper over $\Z_p$, so $X(\Q_p) = \widetilde{\mathcal X}(\Z_p)$ is non-empty if and only if $\widetilde{\mathcal X}_{\F_p}$ has a smooth $\F_p$-point ($\Z_p$ points on a regular model reduce to smooth $\F_p$-points, and conversely, those lift by Hensel's lemma). Away from the nodes $\widetilde{\mathcal X} = \mathcal X$, and over a node $y$ of thickness $t$ lies a chain of $t-1$ copies of $\P^1$, which Frobenius maps to the chain over $\Frob\, y$. This chain has $\F_p$-points only if $\Frob\, y = y$. In that case, if $y$ is split, each $\P^1$ in the chain is defined over $\F_p$ and has $(p+1)-2 \geqslant 1$ smooth $\F_p$-points. If $y$ is non-split, Frobenius reverses the chain, so only the middle $\P^1$, which exists if and only if $t$ is even, is defined over $\F_p$ and consequently has smooth $\F_p$-points.

It remains to translate this to $C_s(\F_p)$, where $C_s = (X_0^D(M)/W')_{\F_p}$. The map $C_s \to \mathcal X_{\F_p}$ induced by $s_1$ is an isomorphism over the smooth locus onto a component (the other one, if $W = W'$, is the image of $C_s$ via $s_2$, with the same supersingular points), and a node is $\F_p$-rational and split if and only if one of its branches, which can be taken to come from $C_s$ via $s_1$, is the image of some $P \in C_s(\F_p)$.
\end{proof}

\begin{remark}
In our implementation of local point checks in \repo{AuxiliaryFiles/existence_of_Qp_points.m}{\texttt{existence{\_}of{\_}Qp{\_}points.m}}, the nodes of the special fiber of $X_0^D(N)$ are the supersingular points of $X_0^D(M)$, i.e., are in correspondence with a basis of the Brandt module in $B_{Dp,\infty}$ (quaternion algebra of discriminant $Dp$). Frobenius permutes these nodes as $w_p$ does (by \cref{rem: comparison of Frob with wp}), and the image of a node $x$ in $X_0^D(N)/W$ is $\Fp$-rational if and only if $\Frob\, x = h\, x$ for some $h \in W$, in which case it is split if and only if this is true for some $h\in W'$.

If $W=W'$ (two component case), every $\Fp$-rational node is split and comes from exactly one point of $C_s(\Fp)$ so it suffices to find an $\Fp$-rational node of thickness at least $2$ and otherwise to compare $\#C_s(\Fp)$ with the number of $\Fp$-rational nodes (in order to detect nodes that are smoothed out on the quotient). If $W\neq W'$ (one component case), a split $\Fp$-rational node has two $\F_p$-rational preimages in $C_s$ and a non-split has none, so we look for an $\Fp$-rational node of even thickness, or a split one of odd thickness at least $3$, and otherwise we compare $C_s(\F_p)$ with twice the number of split $\Fp$-rational nodes.
\end{remark}


\subsection{An interlude: some equations of Atkin--Lehner quotients}

In \cite[Theorem 8.3]{PS25}, the authors determine defining equations for many Atkin--Lehner quotients $X = X_0^D(N)/W$, with $D>1$, such that $X$ has genus $2$ and admits a bielliptic Atkin--Lehner involution. In a smaller number of cases, they determine finitely many candidate equations, exactly one of which must give the correct isomorphism class over $\mathbb{Q}$. These determinations were then improved by Assaf--Hashimoto in \cite{AH26} using the authors' more general implementation of computations of residue fields of CM points (coming from work of Gonz\'{a}lez--Rotger \cite[Appendix]{GR06}). 

Using our implemented reduction type checks, we are able to make the following further improvement on these results. 

\begin{proposition}
\label{prop : bielliptic_improvement}
Each of the following four bielliptic genus $2$ curves is isomorphic over $\mathbb{Q}$ to the hyperelliptic curve with the stated model:
\begin{align*}
    X_0^{14}(45)/\langle w_9, w_{10}, w_{14} \rangle: \quad y^2 &= -2646x^6 + 3807x^4 + 684x^2 + 27, \\
    X_0^{15}(44)/\langle w_4, w_{15}, w_{33} \rangle: \quad y^2 &= 1215x^6 - 2628x^4 + 387x^2 + 18, \\
    X_0^{15}(68)/\langle w_4, w_{15}, w_{17} \rangle: \quad y^2 &= -81x^6 + 1044x^4 - 1989x^2 + 450, \\
    X_0^{34}(15)/\langle w_3, w_{10}, w_{17} \rangle: \quad y^2 &= -432x^6 + 369x^4 + 198x^2 + 9.
\end{align*}
\end{proposition}
\begin{proof}
For each listed curve $X = X_0^D(N)/W$, we have from the work referenced above finitely many candidate equations for $X$ over $\Q$. We find that among these equations, only that listed in the statement of the proposition has the same reduction type as $X$ at all odd primes $p \parallel N$. These computations for the candidate hyperelliptic equations are carried out using Dokchitser's package \cite{redlib}, while the computation for $X$ itself is handled using our implementation of \cref{cor:minimal_regular_AL_quot} in the file \repo{AuxiliaryFiles/dual_graphs_level.m}{dual{\_}graphs{\_}level.m}. 
\end{proof}

\begin{remark}
Working from the final results on bielliptic genus $2$ curves from \cite{PS25}, our reduction type computations at odd primes $p \parallel N$ also recover $70$ of the additional $168$ equations determined in \cite{AH26}. The $4$ equations listed in \cref{prop : bielliptic_improvement} are to our knowledge previously unknown. Our methods also in some cases further narrow the list of potential candidate equations compared to the results of \cite{PS25} and \cite{AH26}, and for this we refer the reader to the file
\repo{Genus2Bielliptics/genus_2_bielliptics_eqn_not_determined_APS.m}{genus{\_}2{\_}bielliptics{\_}eqn{\_}not{\_}determined{\_}APS.m}.
\end{remark}

\begin{example}
By the results of \cite{PS25}, the curve $X \colonequals X_0^{34}(15)/\langle w_3,w_{10},w_{17} \rangle$ is either isomorphic over $\mathbb{Q}$ to the hyperelliptic curve $H_1\colon y^2 = -432x^6 + 369x^4 + 198x^2 + 9$ from \cref{prop : bielliptic_improvement} or to the curve
\[ H_2 : \quad y^2 = \dfrac{27648}{25} x^6 - \frac{7623}{25}x^4 - \frac{2259}{25}x^2 + \frac{576}{25}. \]
Using Dokchitser's package, we compute that the minimal regular models of the reduction types of $H_1$ and $H_2$ at $p=5$ are distinct; that of $H_1$ is $I_{2-2-0}$, hence has dual graph
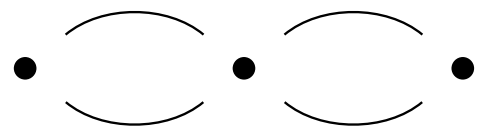
\[\begin{tikzcd}
	\bullet & \bullet & \bullet
	\arrow[curve={height=-12pt}, no head, from=1-2, to=1-1]
	\arrow[curve={height=12pt}, no head, from=1-2, to=1-1]
	\arrow[curve={height=-12pt}, no head, from=1-2, to=1-3]
	\arrow[curve={height=12pt}, no head, from=1-2, to=1-3]
\end{tikzcd}\]
while that of $H_2$ is $I_2-I_2-2$, with the following dual graph:
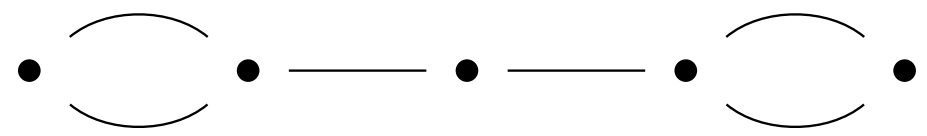
\[\begin{tikzcd}
	\bullet & \bullet & \bullet & \bullet & \bullet
	\arrow[curve={height=-12pt}, no head, from=1-2, to=1-1]
	\arrow[curve={height=12pt}, no head, from=1-2, to=1-1]
	\arrow[no head, from=1-3, to=1-2]
	\arrow[no head, from=1-4, to=1-3]
	\arrow[curve={height=-12pt}, no head, from=1-4, to=1-5]
	\arrow[curve={height=12pt}, no head, from=1-4, to=1-5]
\end{tikzcd}\]
The former matches the reduction type of $X$. 
\end{example}

Similarly, in \cite[Theorem 7.2]{PS25} the authors determine defining equations for many non-elliptic genus $1$ Atkin--Lehner quotients $X = X_0^D(N)/W$. In four of the cases studied, two candidate equations are determined and the authors' methods are not able to determine which is correct. Of these, two have $N=1$, but in the other $2$ cases $N$ is prime and the reduction type at $N$ is sufficient to determine the isomorphism class, giving the following result. 
\begin{proposition}
\label{prop : non-elliptic_improvement}
The following non-elliptic genus $1$ curves are isomorphic over $\mathbb{Q}$ to the curve with the stated model:
\begin{align*}
    X_0^{14}(11)/\langle w_2, w_{11} \rangle:& \quad y^2 = -11x^4 + 144342x^2 - 16964640x - 2075574699, \\
    X_0^{46}(3)/\langle w_2, w_3 \rangle:& \quad y^2 = -3x^4 - 5994x^2 + 59616x - 1384371. 
\end{align*}
\end{proposition}

\begin{proof}
In both cases (at $p=11$ and $p=3$, respectively), the other candidate equation from \cite[Table~8, Theorem~1.4]{PS25} gives a non-reduced special fiber of the minimal regular model, which for this type coincides with \texttt{redlib}'s model $[2]I_4$, (i.e., a fiber with components of multiplicity $>1$). This is not possible by \Cref{rem: genus1}.

The stated equation in each case gives a reduced one ($I_4$), which matches that computed via the moduli interpretation. 
\end{proof}


\section{Isomorphic Atkin--Lehner quotients}\label{section: isomorphisms}     

If two Atkin--Lehner quotients of Shimura curves $X_0^{D_1}(N_1)$ and $X_0^{D_2}(N_2)$ are geometrically isomorphic, then their Jacobians are necessarily isogenous and thus this pair of quotients is referenced in \cref{theorem: main_isogenies_theorem}. In this section, we make progress on determining which of these quotients are isomorphic over $\mathbb{Q}$. While some of our techniques work to show that certain quotients are not \emph{geometrically} isogenous, most of our tools are only sufficient for investigating isomorphism over $\Q$. We thus leave the general question of the existence of isomorphisms over $\overline{\Q}$ to future work. 

\begin{remark}\label{remark: exceptional}
For good reason, modular maps between modular curves are often of particular interest. For example, in \cite[Theorem 4, Section 7]{DHNS26} such maps serve as part of the authors' geometric explanation of rational points on all modular curves, assuming \cite[Conjecture 4.7]{DHNS26}. In this work, however, when we discuss isomorphisms between Atkin--Lehner quotients we mean isomorphisms as algebraic curves with no additional structure assumed. To be clear: if $X/\Q$ and $Y/\Q$ are Atkin--Lehner quotients of, say, modular curves, then we have natural quotient maps $\mathbb{H} \to X(\C)$ and $\mathbb{H} \to Y(\C)$. When we talk about $X$ and $Y$ as algebraic curves in this work we consider them simply as schemes, rather than as weighted projective schemes or as stacks, and we are forgetting about their uniformizations by $\mathbb{H}$; the existence of an isomorphism $\phi : X \to Y$ for us \emph{does not} imply the existence of an automorphism $\alpha$ of $\mathbb{H}$ so that 
\[
\begin{tikzcd}
\mathbb{H} \arrow[r,"\alpha"] \arrow[d] & \mathbb{H} \arrow[d] \\
X(\C) \arrow[r,"\phi"] & Y(\C)
\end{tikzcd}
\]
commutes. In particular, $\phi$ need not preserve the property of being a cuspidal point or of being a CM point. 
\end{remark}

We begin with a lemma handling the most laborious isomorphism check.

\begin{lemma}\label{lemma: annoying_genus_four_iso}
    The curves $X_1 \colonequals X_0(93)/\langle w_{93}\rangle$ and $X_2 \colonequals X_0(186)/\langle w_3, w_{62}\rangle$ are isomorphic over $\Q$.
\end{lemma}    
\begin{proof} We tried confirming this isomorphism of genus $4$ curves in Magma, but the computation did not finish. Hence we present the following direct proof. 
Fix a prime $p$ of good reduction. We compute an explicit $\F_p$-isomorphism $\varphi_p : X_{1,\F_p}\longrightarrow X_{2,\F_p}$.
Since both curves are embedded in $\P^3$, such an isomorphism is induced by a projective linear change of coordinates, hence by a matrix $M_p\in \PGL_4(\F_p)$. After choosing a normalization (fixing the first nonzero entry of $M_p$ to be $1$), we obtain a well-defined representative $M_p$. Repeating this for several primes $p_1,\dots,p_r$ yields a family of matrices $\{M_{p_i}\}$.

We then lift these matrices entrywise to characteristic zero. Viewing $M_{p_i}$ as a $4\times4$ matrix with entries in $\Z/p_i\Z$, we combine the congruences via the Chinese remainder theorem (CRT) to obtain a matrix $M$ modulo $N=\prod \limits_i p_i$ whose reduction modulo each $p_i$ agrees with $M_{p_i}$. This produces an integer matrix $A$ with the same reductions. A direct check shows that $A$ indeed sends the defining equations of $X_1$ to those of $X_2$, and therefore induces a $\Q$-isomorphism $X_1\cong X_2$.

\end{proof}

\begin{proof}[Proof of \cref{theorem: main_isomorphisms_theorem}]
For each pair of quotients that are claimed to be isomorphic in the statement, we have models for both quotients over $\mathbb{Q}$ and we explicitly compute that the curves in the pair are isomorphic over $\mathbb{Q}$. In the case $D_i = 1$ these models can be computed using the Magma \texttt{ModularCurveQuotient} intrinsic, while for $D_i > 1$ these models are found either in work of Guo--Yang \cite{GY} and Padurariu--Schembri \cite{PS23}, namely in the case where the top curve $X_0^{D_i}(N_i)$ is hyperelliptic, or in \cite[Theorem 1.4 (2)]{PS25}, which provides models for certain bielliptic genus $2$ quotients. 

Moreover, for \emph{all} pairs of quotients with isogenous Jacobians for which we have models over $\Q$ for both quotients, we run explicit isomorphism checks before proceeding with further checks. It is worth quickly remarking more on these explicit checks. In many cases, we first check for isomorphism of reductions over $\F_p$ for primes $p$ of good reduction for both quotients. For all cases in which we have equations for both quotients, aside from the single pair dealt with in \cref{lemma: annoying_genus_four_iso}, we are either able to prove non-isomorphism over $\F_p$ for some such $p$ by direct computation in Magma, or else we are able to determine the existence of an isomorphism over $\Q$ directly in Magma. In the case of \cref{lemma: annoying_genus_four_iso}, a direct verification of $\Q$-isomorphism in Magma seemed to be impractical due to excessive memory demands, leading to the alternate argument provided. 

For pairs with isogenous Jacobians for which we lack a model over $\mathbb{Q}$ for at least one quotient, we next run the following filters in an attempt to exclude isomorphism over $\Q$ via local data. 

\begin{enumerate}
    \item \textbf{Local points:} We compare the existence of real points on quotients. Of course, every quotient of a modular curve $X_0(N)$ has real points. For $D_i > 1$, we have that $X_0^{D_i}(N_i)$ has no real points \cite[Theorem 0]{Sh75}, while a result of Ogg \cite[Proposition 1]{Ogg83} determines the existence of real points on a quotient of the form $X_0^{D_i}(N_i)/\langle w_m \rangle$ with $m>1$. That $(X_0^{D_i}(N_i)/W_i)(\R) \neq \emptyset$ is then equivalent to there being some $w_m \in W_i$ with $\left(X_0^{D_i}(N_i)/\langle w_m \rangle\right)(\R) \neq \emptyset$, since $[\C : \R] = 2$ (and the quotient map $X_0^D(N) \rightarrow X_0^D(N)/W$ is defined over $\R$). We also compare the existence of $\Q_p$-points, but, following real point comparisons and the reduction type comparisons described in the next item, such checks did not exclude any additional candidates.
    
    \item \textbf{Reduction types at finite primes:} \cite{PS25} provides an algorithm for computing data of the reduction type of an Atkin--Lehner quotient $X_0^D(N)/W$ (with $N$ squarefree, as in this study) at a prime $p \mid D$. For primes $p \mid N$, we have an analogous implementation, again for arbitrary Atkin--Lehner quotients, coming from \cref{thm: semistable_XD_quotient} and \cref{cor:minimal_regular_AL_quot} in the file \repo{AuxiliaryFiles/dual_graphs_level.m}{dual{\_}graphs{\_}level.m}. If we can find a prime $p$ at which the reduction types disagree, then we will have proven that the quotients are not isomorphic over $\mathbb{Q}$. 
\end{enumerate}

Through these checks, we arrive at the remaining candidate pairs for $\Q$-isomorphism in the referenced tables. 
\end{proof}

\begin{remark}
Code for all computations referenced in the proof of \cref{theorem: main_isomorphisms_theorem} is contained in the files 
\begin{align*}
&\repo{Isomorphisms/SameProductLevel/same_product_isomorphisms_narrowing.m}{\texttt{same{\_}product{\_}isomorphisms{\_}narrowing.m}} \quad \text{and} \\ 
&\repo{Isomorphisms/DifferentProductLevel/different_products_isomorphisms_narrowing.m}{\texttt{different{\_}products{\_}isomorphisms{\_}narrowing.m}}. 
\end{align*}
While we ran these computations separately for the list of isogenous Jacobians of the same levels $D_1N_1 = D_2N_2$ and for the list of quotients of distinct levels $D_1N_1 \neq D_2N_2$, there is no difference in the checks that we apply to exclude the possibility of, or to confirm, isomorphism between these two cases. The true reason for the separation is that we had good reason for separating our computations related to the existence of isogenies (as done in \cref{section: isogenies}), and our lists of pairs of isogenous curves from each of these cases provided a natural separation for our lengthy computations related to isomorphisms. 
\end{remark}

\begin{remark}
As noted in \cite[Remark 6.6]{PS25}, the equations given in Guo--Yang for several pairs of involutions on Shimura curves were swapped, and these swappings carried through to the equations for quotients given in \cite{PS23}. For example, the equations for the involutions $w_2$ and $w_{26}$ on the curve $X_0^{39}(2)$, and of their respective quotients, are swapped in the referenced works, which is important to note in determining the only isomorphism between genus $3$ curves given in the above result. 
\end{remark}

\begin{remark}
There are two further checks for isomorphism we considered, but we do not remark on them in the proof above as they do not give any additional information compared to the checks described therein.

The first check concerns automorphism groups: if we know the (size of) the automorphism group of at least one of the quotients in the pair, then we can compare (the orders of) these groups in a hope to show that they are distinct. Usually the automorphism group just consists of Atkin--Lehner involutions; work of Bars--Gonz\'{a}lez \cite{BG21} determines exactly when $X_0(N)^*$ has a non-trivial automorphism for $N$ squarefree, which would be necessary for a cover $X_0(N)/W$ to have a non-Atkin--Lehner automorphism, and \cite[Section 5]{MPSS25} provides similar checks in the $D>1$ case. Note that this check in fact would in theory exclude the possibility of isomorphism over $\overline{\mathbb{Q}}$, rather than just over $\mathbb{Q}$, as for squarefree $N$ all automorphisms of $X_0^D(N)/W$ are defined over $\Q$ (see \cite[Proposition 1.5]{KR08}). This described check on automorphism groups yields no additional information even prior to our local data considerations; using our available data we could not show via automorphism group information that any pair of curves with isogenous Jacobians over $\Q$ was non-isomorphic.

The second is as follows: if $X$ and $Y$ are isomorphic quotients and $w$ is an Atkin--Lehner involution on $X$, then $Y$ must also have an involution over $\Q$ whose fixed points have the same residue fields over $\Q$ as the fixed points of $w$. For example, consider the following candidate pair $(X,Y)$ of genus $2$ curves for isomorphism over $\mathbb{Q}$:
\[ X := X_0(165)^* \quad \text{ and } \quad Y := X_0^{55}(3)/ \langle w_3, w_{55} \rangle . \]
Using code of Assaf--Hashimoto from \cite{AH26}, specifically that implementing work of Gonz\'{a}lez--Rotger \cite[Appendix]{GR06} on residue fields of CM points on Atkin--Lehner quotients, we determine that the residue fields of the fixed points of the non-trivial involution on $Y$ are as follows: there are two closed points with residue field $\mathbb{Q}$, and the remaining two closed points are quadratic with residue fields $\Q(\sqrt{5})$ and $\Q(\sqrt{-11})$. The modular curve quotient $X$ is isomorphic over $\mathbb{Q}$ to the curve with hyperelliptic model
\[ y^2 = f(x) = x^6 - 4x^4 + 10x^3 - 16x^2 + 9, \]
and via the factorization of $f$ over $\mathbb{Q}$ we see that the fixed points of the non-trivial involution on $Y$ have the same residue fields as that of $X$. 
\end{remark}

In two of the remaining cases from \cref{table: remaining_isomorphism_candidates_same_product}, while we are not able to determine isomorphism we are able to narrow down the possibilities for the isomorphism classes of the relevant quotients.

\section{Zukunftsmusik}\label{section: Zukunftsmusik}

\subsection{Non-squarefree level}\label{section: non-squarefree}

In this work, we have restricted our attention to squarefree level. Note, though, that the proofs of the results we used in \cref{section: isogenies} to reduce our consideration of isogenies between Atkin--Lehner quotients to finitely many cases -- namely, \cref{prop: AL_restrictions}, \cref{prop: AL_restrictions_same_curve}, \cref{prop: dim_zero_prop_shimura}, \cref{prop: general_dim_zero_prop_shimura}, and \cref{level_change_prop} -- do not rely on the levels being squarefree. Rather, the statements of these results involve hypotheses or conclusions related to zero-dimensional spaces of newforms with prescribed signed patterns, and so we only used them in the squarefree case as here we have the complete classification of this phenomenon from \cite{PPV26}. Ergo, an extension of the work of Padurariu--Park--Voight to classify all such spaces of dimension $0$ in general level would help lead to an extension of \cref{theorem: main_isogenies_theorem} to general level. 

It is important to note, though, that unlike the situation for squarefree levels, general level \emph{does not} admit finitude of isogenies between Atkin--Lehner quotients. In particular, Furumoto--Hasegawa observed infinite families of \emph{isomorphic} Atkin--Lehner quotients of the modular curves $X_0(N)$.  
\begin{proposition}\cite[Propositions 4 \& 5]{FH99}\label{FH_isos}
\begin{enumerate}
    \item If $2 \nmid M$ and 
\[ W = \langle w_4, w_{m_1}, \ldots, w_{m_r} \rangle \]
is a subgroup of the Atkin--Lehner group of $X_0(4M)$ with $m_i \parallel M$ for all $1 \leq i \leq r$, then there is an isomorphism over $\mathbb{Q}$
\[ X_0(4M)/W \cong X_0(2M)/\langle w_{m_1}, \ldots, w_{m_r} \rangle. \]

\item Suppose that $3 \nmid M$ and let
\[ W = \langle w_{m_1}, \ldots, w_{m_r} \rangle \]
be a subgroup of the Atkin--Lehner group of $X_0(9M)$ with $m_i \parallel 9M$ for each $1 \leq i \leq r$. Letting
\[ \epsilon(n) = \begin{cases} 0 &\text{ if } n \equiv 1 \pmod{3} \text{ or if } 9 \parallel n \text{ and } n/9 \equiv 1 \pmod{3} \\
1 & \text{ otherwise} \end{cases}, \]
we have an isomorphism over $K = \mathbb{Q}(\sqrt{-3})$
\[ \left(X_0(9M)/W\right) \otimes_{\Spec \Q} \Spec(\Q(\sqrt{-3})) \cong X_0(9M)/\left\langle \left\{W_9^{\epsilon(N_i)} W_{N_i}\right\}_{i=1}^{r} \right\rangle \otimes_{\Spec \Q} \Spec(\Q(\sqrt{3})) . \]
\end{enumerate}
\end{proposition}
Assaf--Hashimoto, using the fact that relevant congruence subgroup arithmetic is the same in the $D>1$ case, extend this observation to analogous quotients of Shimura curves \cite[Proposition 6.1.1]{AH26}. For clarity of what we have discussed above, we can reframe the first part of \cref{FH_isos} as giving a statement about spaces of newforms with prescribed sign patterns. 
\begin{proposition}

    Let $M$ be an odd integer, and let $\overline{\epsilon}$ be a sign pattern in level $4M$ with sign $\epsilon_1$ at $2$. If 
    \[ \dim\left(S_2(4M)^{\overline{\epsilon}}\right) > 0, \]
    then $\epsilon_1 = -1$.
\end{proposition}
\begin{proof}
By \cref{FH_isos}, have an isomorphism
    \[ X_0(4M)/\langle w_4 \rangle \cong X_0(2M). \]
    Therefore, every isogeny factor of $J = \Jac(X_0(4M)/\langle w_4 \rangle)$ is associated to a newform of level dividing $2M$. Hence, any newform $f$ of level $4M$ has sign $-1$ at $2$. 
\end{proof}
The second part of \cref{FH_isos} does not, in general, lead to a result on newforms over $\Q$ because the isomorphisms described there are not defined over $\mathbb{Q}$. For example, $X_0(45)/\langle w_5 \rangle$ and $X_0(45)/\langle w_{45} \rangle$ are distinct twists of the same elliptic curve over $\mathbb{Q}$ which become isomorphic upon base change to $\Q(\sqrt{-3})$. 

At best, one could then hope for a result which parameterizes all families of isogenies (over $\Q$ or geometrically) between Jacobians of Atkin--Lehner quotients outside of an explicit, finite list of ``sporadic'' isogenies in the general level case. We leave this pursuit to future work, but natural (closely related) questions to begin with could be as follows. 

\begin{question}
Are there only finitely many isomorphisms (or more generally, isogenies) over $\mathbb{Q}$ between distinct Atkin--Lehner quotients of modular curves of the form $X_0(N_1)$ and $X_0(N_2)$ with $N_1N_2$ odd? 
\end{question}

\begin{question}
Are there finitely many pairs $(N,\overline{\epsilon})$, with $N$ odd and $\overline{\epsilon}$ a sign pattern in level $N$, for which
\[ \dim\left(S_2(N)^{\overline{\epsilon}}\right) = 0 ? \]
\end{question}

\subsection{General congruence subgroups}

In this work we have focused on Atkin--Lehner quotients of the Shimura curves $X_0^D(N)$. Of course, one can more ambitiously ask about the existence of isomorphisms (or isogenies of Jacobians) between quotients of more general modular or Shimura curves attached to specified subgroups $H \leq \GL_2(\Z/N\Z)$. Perhaps the best-known example is the isomorphism $X_{\mathrm{ns}}^+(13) \cong_\Q X_{\mathrm{s}}^+(13)$ between the modular curves
attached to the normalizers of the non-split and split Cartan subgroups of level~$13$, proved by Baran~\cite{Baran14}. This isomorphism is not modular. Further examples appear sporadically in the literature. For
instance, the genus-$2$ modular curve with LMFDB label \href{https://www.lmfdb.org/ModularCurve/Q/18.108.2.a.1/}{\texttt{18.108.2.a.1}} is isomorphic over $\Q$ to $X_1(18)$, although the isomorphism does not respect the maps to the $j$-line; see \cite[Proposition~3.5 and Appendix~A]{AG26}.

As we have done in this work, it would be reasonable to restrict to quotients of Shimura curves $X_\Gamma$ attached to a congruence subgroup $\Gamma$ by subgroups of $\Aut(X)$ consisting of \emph{modular automorphisms}. That is, one may restrict to the subgroup
\[ B(\Gamma) := \mathcal{N}(\Gamma)/\Gamma \leq \Aut_{\mathbb{C}}\left(X_\Gamma\right), \] 
where $\mathcal{N}(\Gamma)$ denotes the normalizer of $\Gamma$ in $\PSL_2(\Z/N\Z)$. The expectation is that $\Aut\left(X_\Gamma\right) = B(\Gamma)$ for $X_\Gamma$ of large enough genus (see \cite[Conjecture 1.1]{KR08}), so this should not be viewed as a significant restriction. Indeed, this expected statement has been proven for certain classical families of modular curves, including $X(N)$, $X_1(N)$ for $N$ squarefree, and $X_0(N)$ (we refer to \cite[Theorem 1.4]{KR08} where the relevant references are collected.\footnote{Note that the reference for $X_1(N)$ is to work of Momose that has unfortunately never appeared in the literature, while the extension of the result of Serre for $X(p)$, for prime $p$, to composite level $N$ is contained in the thesis of Kontogeorgis (written in Greek) but can also be found in \cite[Section 3]{BKX13}.}) Dose--Lido--Mercuri have recently proven the claim for general modular curves $X_H$, with $-I \in H$ and with $H$ having surjective determinant, of prime level \cite[Theorem 3.13]{DLM26}. In the $D>1$ case, one lacks cusps and thus immediately has the claim in the absence of elliptic points \cite[Proposition 1.2]{KR08}, which occurs for example for the analogs $X_1^D(N)$ and $X^D(N)$ of $X_1(N)$ and $X(N)$. For the curves $X_0^D(N)$, we are not aware of a general result even in the squarefree case, but we refer to \cite{KR08} for related progress. 

Whether working with general congruence subgroups $\Gamma$ or with $\Gamma$ restricted to some specified family, as done in this work for the subgroups $\Gamma_0^D(N)$ corresponding to the Shimura curves $X_0^D(N)$, it seems that the main work to be done would be as follows:
\begin{itemize}
    \item fully work out the action of $B(\Gamma)$ on the isogeny factors of $\Jac\left(X_\Gamma\right)$, and
    \item provide an analog of the results of \cite{PPV26} (see \cref{section: PPV}) for the congruence subgroups in question, i.e., classify exactly those triples $(\Gamma,N,\overline{\epsilon})$, with $\overline{\epsilon}$ a sign pattern in level $N$, such that $\dim \left(S_2(\Gamma,N)^{\overline{\epsilon}}\right) = 0$.
    
\end{itemize}
Following such work, our line of argumentation in \cref{section: isogenies} should be sufficient to prove explicit finitude of isogenies between Jacobians of curves ${X_\Gamma}/W$ with $W \leq B(\Gamma)$ a subgroup of modular automorphisms. 

We expect that the modular curves $X_1(N)$ of squarefree level $N$ would be a reasonable testing ground for extending our results and those of \cite{PPV26}. One can find the structure of $B(\Gamma_1(N))$ worked out in \cite{Lang01} (and also independently in \cite{KK00}); it is an extension of the subgroup of Atkin--Lehner involutions by $\left(\Z/N\Z\right)^\times/\{\pm 1\}$ (see \cite[Theorem 1, Corollary 3]{KK00}). After gaining an explicit understanding of the action of this group on the Jacobian, a first step towards the second item above in this case, which would itself be of interest, may be as follows.
\begin{question}
For which pairs $(N,p)$, with $N$ squarefree and $p \mid N$, does the quotient $X_1(N)/B(\Gamma_1(N))$ have good reduction at $p$? 
\end{question}

\section{Tables}

\rowcolors{2}{white}{gray!20}
\begin{longtable}{|c|c|c|c|c|}
\caption{Isogenous Jacobians of distinct Atkin--Lehner quotients $X_0^D(N)/W$ of a fixed Shimura curve $X_0^D(N)$, as discussed in \cref{subsection: isogenies_same_curve}. In each row, the final column lists subgroups $W \leq W_0(D,N)$ of automorphisms of the specified curve $X_0^D(N)$ such that the corresponding quotients $X_0^D(N)/W$ are of genus $g$ with Jacobians in the same isogeny class.}\label{table: isogenies_same_curve} \\ \hline
\rowcolor{gray!50}
$D\cdot N$ & $D$ & $N$ & $g$ & $\textnormal{Subgroups } W \leq W_0(D,N)$ \\ \hline
$ 30 $ & 
$ 1 $ & $ 30 $ & $ 2 $ & 
$\langle 
w_{ 2 } \rangle,$ 
$\langle 
w_{ 10 } \rangle$ \\ \hline
& 
$ 15 $ & $ 2 $ & $ 2 $ & 
$\langle 
w_{ 2 } \rangle,$ 
$\langle 
w_{ 6 } \rangle$ \\ \hline
$ 42 $ & 
$ 1 $ & $ 42 $ & $ 2 $ & 
$\langle 
w_{ 3 } \rangle,$ 
$\langle 
w_{ 21 } \rangle$ \\ \hline
& & & & 
$\langle 
w_{ 6 } \rangle,$ 
$\langle 
w_{ 42 } \rangle$ \\ \hline
& 
$ 14 $ & $ 3 $ & $ 2 $ & 
$\langle 
w_{ 3 } \rangle,$ 
$\langle 
w_{ 6 } \rangle$ \\ \hline
& 
$ 21 $ & $ 2 $ & $ 2 $ & 
$\langle 
w_{ 6 } \rangle,$ 
$\langle 
w_{ 14 } \rangle$ \\
\hline
\end{longtable}

\rowcolors{2}{white}{gray!20}
\begin{longtable}{|c|c|>{\raggedright\arraybackslash}p{30em}|}
\caption{Isogenous Jacobians of quotients $X_0^D(N)/W$ of distinct Shimura curves $X_0^D(N)$ of the same product level $L=DN$ of fixed genus $g > 1$, as discussed in \cref{subsection: isogenies_same_level}. Each row gives quotients in a fixed isogeny class in the designated level $L$ and genus $g$, with each quotient listed as a tuple $(D,N,W)$.}\label{table: isogenies_same_product} \\ \hline
\rowcolor{gray!50}
$L$ & $g$ & $\textnormal{Isogenous curves}$ \\ \hline
$ 26 $ &  $ 2 $ & 
$( 1 , 26 , \langle 
\textnormal{id} \rangle ), $
$( 26 , 1 , \langle 
\textnormal{id} \rangle )$ \\ \hline
$ 30 $ &  $ 2 $ & 
$( 1 , 30 , \langle 
w_{ 2 } \rangle ), $ 
$( 1 , 30 , \langle 
w_{ 10 } \rangle ), $ 
$( 15 , 2 , \langle 
w_{ 2 } \rangle ), $ 
$( 15 , 2 , \langle 
w_{ 6 } \rangle )$ \\ \hline
& & 
$( 1 , 30 , \langle 
w_{ 3 } \rangle ), $ 
$( 15 , 2 , \langle 
w_{ 5 } \rangle )$ \\ \hline
&  $ 3 $ & 
$( 1 , 30 , \langle 
\textnormal{id} \rangle ), $
$( 15 , 2 , \langle 
\textnormal{id} \rangle )$ \\ \hline
$ 35 $ &  $ 2 $ & 
$( 1 , 35 , \langle 
w_{ 7 } \rangle ), $ 
$( 35 , 1 , \langle 
w_{ 5 } \rangle )$ \\ \hline
&  $ 3 $ & 
$( 1 , 35 , \langle 
\textnormal{id} \rangle ), $
$( 35 , 1 , \langle 
\textnormal{id} \rangle )$ \\
\hline
\end{longtable}

\setlength{\LTcapwidth}{0.95\textwidth}
\begingroup
\begin{longtable}{|c|>{\raggedright\arraybackslash}p{40em}|}
\caption{Isogenous Jacobians of Shimura curve quotients $X_0^D(N)/W$ of genus $g > 1$ of distinct product levels, as discussed in \cref{subsection: isogenies_different_level}. Each row gives members of a fixed isogeny class of quotients, where each quotient is listed as a triple \( (D, N, W)\). Note, however, that as we record all such members of a fixed isogeny class per row, some rows also contain pairs with same product levels appearing in prior tables.}\label{tab:isog_shimura_different_levels}\\
\hline
$g$ & Isogeny class \\
\hline
\endfirsthead
\hline
$g$ & Isogeny class \\
\hline
\endhead
$2$ & $(1,23,\langle \textnormal{id}\rangle)$, $(1,46,\langle w_{46}\rangle)$, $(1,69,\langle w_{69}\rangle)$ \\ \hline
 & $(1,26,\langle \textnormal{id}\rangle)$, $(26,1,\langle \textnormal{id}\rangle)$, $(1,78,\langle w_{39}\rangle)$, $(26,3,\langle w_{78}\rangle)$, $(26,3,\langle w_{6}\rangle)$, $(26,5,\langle w_{65}\rangle)$, $(26,7,\langle w_{182}\rangle)$, $(26,15,\langle w_{6}, w_{65}\rangle)$ \\ \hline
 & $(1,29,\langle \textnormal{id}\rangle)$, $(1,87,\langle w_{87}\rangle)$ \\ \hline
 & $(1,30,\langle w_{10}\rangle)$, $(1,30,\langle w_{2}\rangle)$, $(15,2,\langle w_{2}\rangle)$, $(15,2,\langle w_{6}\rangle)$, $(15,14,\langle w_{6}, w_{7}\rangle)$, $(15,22,\langle w_{6}, w_{55}\rangle)$, $(15,34,\langle w_{2}, w_{255}\rangle)$ \\ \hline
 & $(1,31,\langle \textnormal{id}\rangle)$, $(1,62,\langle w_{62}\rangle)$ \\ \hline
 & $(1,33,\langle w_{3}\rangle)$, $(1,66,\langle w_{3}, w_{22}\rangle)$, $(1,66,\langle w_{2}, w_{3}\rangle)$ \\ \hline
 & $(1,35,\langle w_{7}\rangle)$, $(35,1,\langle w_{5}\rangle)$, $(1,70,\langle w_{2}, w_{7}\rangle)$, $(1,70,\langle w_{7}, w_{10}\rangle)$, $(35,2,\langle w_{5}, w_{14}\rangle)$, $(35,2,\langle w_{2}, w_{5}\rangle)$, $(35,3,\langle w_{3}, w_{5}\rangle)$, $(35,3,\langle w_{5}, w_{21}\rangle)$, $(35,6,\langle w_{3}, w_{5}, w_{14}\rangle)$, $(35,6,\langle w_{2}, w_{3}, w_{5}\rangle)$, $(35,6,\langle w_{2}, w_{5}, w_{21}\rangle)$ \\ \hline
 & $(1,37,\langle \textnormal{id}\rangle)$, $(1,74,\langle w_{74}\rangle)$, $(1,111,\langle w_{111}\rangle)$ \\
\hline
\end{longtable}
\endgroup

\rowcolors{2}{white}{gray!20}
\begin{longtable}{|c|c|c|c|c|c|c|}
\caption{All $68$ pairs of Shimura curves $X \coloneqq X_0^{D_1}(N_1)/W_1$ and $Y \coloneqq X_0^{D_2}(N_2)/W_2$, with $D_1N_1 = D_2N_2$ and with $N_1$ and $N_2$ squarefree, for which we remain unsure of whether $X$ and $Y$ are isomorphic over $\mathbb{Q}$ following \Cref{theorem: main_isomorphisms_theorem}. For each pair we also give the genus $g$ of $X$ and $Y$.}\label{table: remaining_isomorphism_candidates_same_product} \\ \hline
\rowcolor{gray!50}
$D_1$ & $N_1$ & $W_1$ & $D_2$ & $N_2$ & $W_2$ & $g$ \\ \hline
$ 1 $ & $ 114 $ & $ \langle
w_{ 2 }, 
w_{ 19 } \rangle $ & 
$ 57 $ & $ 2 $ & $\langle
w_{ 114 } \rangle $ & $ 
3 $ \\ \hline
$ 1 $ & $ 115 $ & $ \langle
w_{ 5 }, 
w_{ 23 } \rangle $ & 
$ 115 $ & $ 1 $ & $\langle
w_{ 115 } \rangle $ & $ 
2 $ \\ \hline
$ 1 $ & $ 154 $ & $ \langle
w_{ 2 }, 
w_{ 7 }, 
w_{ 11 } \rangle $ & 
$ 77 $ & $ 2 $ & $\langle
w_{ 2 }, 
w_{ 77 } \rangle $ & $ 
2 $ \\ \hline
$ 1 $ & $ 161 $ & $ \langle
w_{ 7 }, 
w_{ 23 } \rangle $ & 
$ 161 $ & $ 1 $ & $\langle
w_{ 161 } \rangle $ & $ 
2 $ \\ \hline
$ 1 $ & $ 165 $ & $ \langle
w_{ 3 }, 
w_{ 5 }, 
w_{ 11 } \rangle $ & 
$ 55 $ & $ 3 $ & $\langle
w_{ 3 }, 
w_{ 55 } \rangle $ & $ 
2 $ \\ \hline
$ 1 $ & $ 186 $ & $ \langle
w_{ 2 }, 
w_{ 3 }, 
w_{ 31 } \rangle $ & 
$ 93 $ & $ 2 $ & $\langle
w_{ 2 }, 
w_{ 93 } \rangle $ & $ 
2 $ \\
\hline
\end{longtable}

\rowcolors{2}{white}{gray!20}
\begin{longtable}{|c|c|c|c|c|c|c|}
\caption{All $83$ pairs of quotients $X \coloneqq X_0^{D_1}(N_1)/W_1$ and $Y \coloneqq X_0^{D_2}(N_2)/W_2$, with $D_1N_1 \neq D_2N_2$ and with $D_1N_1$ and $D_2N_2$ squarefree, of genus $g \geq 2$, for which we remain unsure of whether $X$ and $Y$ are isomorphic over $\mathbb{Q}$ following \Cref{theorem: main_isomorphisms_theorem}. For each pair we also give the genus $g$ of $X$ and $Y$.}\label{table: remaining_isomorphism_candidates_different_products} \\ \hline
\rowcolor{gray!50}
$D_1$ & $N_1$ & $W_1$ & $D_2$ & $N_2$ & $W_2$ & $g$ \\ \hline
$ 6 $ & $ 65 $ & $ \langle w_{2}, w_{3}, w_{13} \rangle $ & 
$ 390 $ & $ 7 $ & $\langle w_{2}, w_{3}, w_{5}, w_{91} \rangle$ & $ 
2 $ \\ \hline
$ 6 $ & $ 91 $ & $ \langle w_{2}, w_{3}, w_{13} \rangle $ & 
$ 546 $ & $ 5 $ & $\langle w_{2}, w_{3}, w_{7}, w_{65} \rangle$ & $ 
2 $ \\ \hline
$ 6 $ & $ 115 $ & $ \langle w_{2}, w_{3}, w_{115} \rangle $ & 
$ 6 $ & $ 161 $ & $\langle w_{2}, w_{3}, w_{161} \rangle$ & $ 
2 $ \\ \hline
$ 10 $ & $ 57 $ & $ \langle w_{2}, w_{5}, w_{19} \rangle $ & 
$ 570 $ & $ 7 $ & $\langle w_{2}, w_{3}, w_{5}, w_{133} \rangle$ & $ 
2 $ \\ \hline
$ 10 $ & $ 77 $ & $ \langle w_{2}, w_{35}, w_{55} \rangle $ & 
$ 770 $ & $ 3 $ & $\langle w_{2}, w_{7}, w_{11}, w_{15} \rangle$ & $ 
2 $ \\ \hline
$ 10 $ & $ 77 $ & $ \langle w_{5}, w_{14}, w_{22} \rangle $ & 
$ 770 $ & $ 3 $ & $\langle w_{5}, w_{6}, w_{7}, w_{11} \rangle$ & $ 
2 $ \\
\hline
\end{longtable}

\bibliographystyle{amsalpha-nodash}
\bibliography{biblio}
\end{document}